%% file: ODMmor.tex
\documentclass[moor]{informs1hack}

\usepackage{bbm, bm}
\usepackage{amsmath, amssymb}
\usepackage{amsfonts}
\usepackage{tabularx}
\usepackage{multicol}
\usepackage{url}
\usepackage{graphicx}
\usepackage{pict2e,color}
\usepackage{dsfont}
\usepackage{bm}
\usepackage{breqn}
\usepackage{colonequals}
\DeclareMathOperator{\Vol}{Vol}

\usepackage{natbib}
 \NatBibNumeric
 \def\bibfont{\small}%
 \def\bibsep{\smallskipamount}%
 \def\bibhang{24pt}%
 \bibpunct[, ]{[}{]}{,}{n}{}{,}%

\usepackage{hyperref}

\TheoremsNumberedBySection  

\EquationsNumberedThrough    

\newcommand{\Po}{\mathcal{P}}
\newcommand{\Q}{\mathcal{Q}}
\newcommand{\T}{\mathcal{T}}

\newcommand{\R}{\mathbb{R}}

\allowdisplaybreaks

\begin{document}



 \RUNAUTHOR{Speakman and Lee}

\RUNTITLE{Quantifying Double McCormick}


\TITLE{Quantifying Double McCormick}



\ARTICLEAUTHORS{%
\AUTHOR{Emily Speakman and Jon Lee}
\AFF{University of Michigan\\ \EMAIL{$\{$eespeakm,~jonxlee$\}$@umich.edu}}
} 

%

%
%
%

%

\ABSTRACT{%
Abstract: When using the standard McCormick inequalities twice to convexify trilinear monomials,
as is often the practice in modeling and software,
there is a choice of which variables to group first. For the important case in which the domain is a nonnegative box,
we calculate the volume of the resulting relaxation, as a function of the bounds defining the box.
In this manner, we precisely quantify the strength of the different possible relaxations defined by all three groupings, in addition to the trilinear hull itself. As a by product, we characterize the best double-McCormick relaxation.

We wish to emphasize that, in the context of spatial branch-and-bound for factorable formulations, our results do not only apply to variables in the input formulation.
Our results apply to monomials that involve auxiliary variables as well.
So, our results apply to the product of any three (possibly complicated) expressions in a formulation.

}

\MSCCLASS{90C26}
\ORMSCLASS{Primary: global optimization; secondary: mixed-integer nonlinear optimization}
\AREAOFREVIEW{Discrete Optimization}

\maketitle


\input{intro.tex}

\input{double.tex}

\input{alternatives.tex}

\input{theorems.tex}

\input{theorem_H.tex}

\input{theorem_A.tex}

\input{theorem_B.tex}

\input{theorem_C.tex}

\input{future.tex}

\input{Appendix.tex}

\section*{Acknowledgments.} The authors gratefully acknowledge partial support from  NSF grant CMMI--1160915, ONR grant N00014-14-1-0315, and a Rackham Summer Award.

\bibliographystyle{ormsv080}
\bibliography{biblioMC} 

\end{document}

%% file: intro.tex
\section{Introduction.}
\label{Trilinear Intro}
Spatial branch-and-bound (sBB) (see \cite{Adjiman98,Ryoo96,Smith99}) for so-called factorable  math\-ematical-optimization formulations (see \cite{McCormick76}) is the workhorse general-purpose algorithm in the area of global optimization. It
works by using additional variables to reformulate every function of the formulation as a (labeled) directed acyclic graph (DAG).
Root nodes can be very complicated functions, and leaves are variables
that appear in the input formulation,
each labeled with its interval domain. Intermediate nodes are labeled with auxiliary variables
together with operators from a small dictionary of basic functions of few (often one or two) variables.
Also, we have a method for convexifying the graph of each dictionary function.
sBB algorithms work by composing convex relaxations of the dictionary functions, according to the DAG,
to get relaxations of the root functions. Bounds on the leaves propagate
to other nodes and conversely. Branching (subdividing the domain interval of a variable)
creates subproblems, which are treated recursively. Objective bounds for subproblems are appropriately combined
to achieve a global-optimization algorithm.

Much of the research on sBB has focused on developing tight convexifications for basic functions of
few variables (many references can be found in \cite{CafieriLeeLiberti10}).
Other research has focused on how bounds can be efficiently propagated
and how branching can be judiciously be carried out (see \cite{Belotti09}, for example).
From the viewpoint of
good convexifications, much less attention has been paid to how the DAGs are created, but this can have a strong impact on the quality of the resulting convex relaxation of the input formulation; see \cite{LibertiCafieri1, LibertiCafieri2,Schichl1,Schichl2} for some key papers with other viewpoints concerning constructing DAGs.

For basic multilinear monomials $f(x_1,\ldots,x_n):=x_1\cdots x_n$, with $x_i\in [a_i,b_i]$,
there is already a lot of flexibility which can have a significant impact on the
overall convexification of the graph of $f(x_1,\ldots,x_n):=x_1\cdots x_n$ on the
box domain $[a_1,b_1]\times \cdots \times [a_n,b_n]$. For $n=2$, we have the classic McCormick inequalities (see \cite{McCormick76}),
which is simply the tetrahedron that is the convex hull of the points $(f,x_1,x_2):=(a_1a_2,a_1,a_2),
(a_1b_2,a_1,b_2), (b_1a_2,b_1,a_2), (b_1b_2,b_1,b_2)$. The inequalities can be derived from the four inequalities
\[
\begin{array}{cc}
  (x_1-a_1)(x_2-a_2)  \geq  0, & (x_1-a_1)(b_2-x_2) \geq  0, \\
  (b_1-x_1)(x_2-a_2)  \geq  0, & (b_1-x_1)(b_2-x_2)  \geq  0,
\end{array}
\]
by multiplying out and then replacing all occurrences of $x_1x_2$ by the variable $f$.

For general $n$, there are $2^n$ points to consider (i.e., all choices of
each variable at a bound), and the inequality descriptions in the space of
$(f,x_1,\ldots,x_n)\in\mathbb{R}^{n+1}$ get rather complicated (even for $n=3$; see \cite{Meyer04a,Meyer04b}).
It is frequent practice,
both in modeling and software, to repeatedly use the McCormick inequalities when $n>2$. Already the trilinear case, $n=3$, is an interesting one for analysis.
Here, we have three choices, which can be thought of as
$f=(x_1x_2)x_3$, $f=(x_1x_3)x_2$ and $f=(x_2x_3)x_1$.
Because the domain of each variable is its own interval
$[a_i,b_i]$, the grouping can affect the quality of the convexification.
In what follows, we analytically quantify
the quality of these different convexification possibilities,
in addition to the trilinear hull itself.

\emph{Our results are not just relevant to trilinear monomials in formulations.}
With the sBB approach for factorable formulations, our results are relevant
whenever three quantities are multiplied. That is, as an expression DAG
is created and auxiliary variables are introduced, a trilinear monomial
will arise whenever three quantities (which can be complicated
functions themselves) are multiplied.

In what follows, we use $(n+1)$-dimensional volume to compare different natural convexifications of graphs of
functions of $n$ variables on the box domain $[a_1,b_1]\times \cdots \times [a_n,b_n]$.
We present a complete analytic analysis of the case of $n=3$, for all choices of $0\leq a_i<b_i$. It is perhaps surprising that this can be carried out,
and probably less surprising that the analysis is quite complicated.

Volume as a measure for comparing relaxations was first proposed in
\cite{Lee94}; also see \cite{KLS1997} and \cite{Steingr}. In fact, the
\emph{practical} use of volume as a measure for comparing relaxations
in the context of nonlinear mixed-integer optimization,
foreshadowed by \cite{Lee94}, was later validated computationally
for a nonlinear version
of the uncapactitated facility-location problem (see \cite{Lee_2007}).
Specifically, using volume calculations,  a main mathematical result of \cite{Lee94} is that weak formulations
of facility-location problems are very close to
strong formulations when the number of  facilities is small compared to the number of customers.
Then \cite{Lee_2007} showed that in this scenario, with a convex objective function,
the weak formulation  computationally out performs the strong formulation in the context of branch-and-bound.
The emphasis in \cite{Lee94,KLS1997,Steingr}
 was not on sBB nor
on low-dimensional functions.
Because those results pertained to varying dimension and
related asymptotics, exactly how volumes are compared and scaled was important (in particular, see \cite{Lee94}
which defines the ``idealized radial distance'').
Because we now focus on low-dimensional polytopes, the exact manner of comparison and scaling is much less relevant. Using volume as a measure
corresponds to a uniform distribution of the optimal solution
across a relaxation. This is
justified in the context of \emph{nonlinear} optimization
if we want a measure that is robust across all formulations.
One can well find situations where the volume measure is misleading.
It would not make sense for evaluating polyhedral relaxations
of the integer points in a polytope, if we were only concerned with \emph{linear}
objectives --- in such a case, solutions are concentrated on the
boundary and there are better measures available (see \cite{Lee94}).
But if we are interested in a mathematically-tractable measure that
robustly makes sense in the
context of global optimization, volume is quite natural.

Motivated by two well-studied applications (the \emph{Molecular Distance
Geometry Problem} and the \emph{Hartree-Fock Problem}),
\cite{CafieriLeeLiberti10} first proposed volume in the context of sBB
and monomials, but they leapfrogged to the case of $n=4$ and took a mostly experimental approach. They  demonstrated that there
can be a significant difference in performance depending on grouping, and
they offered some guidance based on computational experiments.  At the time of that work, it appeared that developing precise formulae for volumes relevant to
repeated McCormick was not tractable. With our present work on $n=3$,
it now seems possible that the case of $n=4$ could be carried out (see \S\ref{sec:future} for an idea concerning how our results for $n=3$
could already be applied in practice to the $n=4$ case).

There has been considerable research on multilinear monomials and generalizations in the context of global optimization,
notably \cite{Rikun97,Luedtke12,BKS2015,RyooSahi,JMW08,EdgeConcave_MeyerFloudas}. Our work adds to that literature.

In \S\ref{sec:double}, we define the polytopes that we work with.
In \S\ref{sec:alternatives}, we discuss the various alternatives for working with triple products.
In \S\ref{sec:theorems}, we present our main results and their consequences.
In \S\S\ref{Three Continuous Convex Hull}--\ref{linC}, we present our proofs.
In \S\ref{sec:future}, we describe future directions for investigation.
\S\ref{app}, an appendix, contains technical lemmas and calculations.

%% file: double.tex
\section{Double McCormick.}\label{sec:double}
\label{Double McCormick Inequalities}

When using the double-McCormick technique to convexify trilinear monomials, a modeling/algorithmic choice is involved: we must choose to which pair of variables we apply the first iteration of McCormick. For the variables $x_i \in [a_i,b_i]$, $i=1,2,3$, let $\mathcal{O}_i \colonequals a_i(b_jb_k) + b_i(a_ja_k)$.
Then we can label the variables such that $\mathcal{O}_1 \leq \mathcal{O}_2 \leq \mathcal{O}_3$.
In this manner, we can assume that
\[\label{Omega}
a_1b_2b_3+b_1a_2a_3 \leq b_1a_2b_3 + a_1b_2a_3 \leq b_1b_2a_3 + a_1a_2b_3. \tag{$\Omega$}
\]

Given the trilinear monomial $f:=x_1x_2x_3$, there are three choices of convexifications depending on the bilinear sub-monomial we  convexify first. We could first group $x_1$ and $x_2$ and convexify $w=x_1x_2$; after this, we are left with the monomial $f=wx_3$ which we can also convexify using McCormick. Alternatively, we could first  group variables $x_1$ and $x_3$, or variables $x_2$ and $x_3$.

\subsection{Convexification.}
\label{lin}

To see how to perform these convexifications in general, we show the double-McCormick convexification that first groups the variables $x_i$ and $x_j$.  Therefore we have $f=x_ix_jx_k$ and we let $w_{ij}=x_ix_j$ so $f=w_{ij}x_k$.

Convexifying $w_{ij}=x_ix_j$ we obtain the inequalities:
\begin{align*}
w_{ij}-a_jx_i-a_ix_j+a_ia_j \geq 0, \\
-w_{ij} + b_jx_i + a_ix_j - a_ib_j \geq 0, \\
-w_{ij} +a_jx_i + b_ix_j - b_ia_j \geq 0, \\
w_{ij} - b_jx_i - b_ix_j + b_ib_j \geq 0.
\end{align*}
Convexifying $f=w_{ij}x_k$ we obtain the inequalities:
\begin{align*}
f- a_kw_{ij} - a_ia_jx_k + a_ia_ja_k \geq 0, \\
-f +b_kw_{ij} + a_ia_jx_k - a_ia_jb_k \geq 0, \\
-f + a_kw_{ij} + b_ib_jx_k - b_ib_ja_k \geq 0, \\
f - b_kw_{ij} - b_ib_jx_k + b_ib_jb_k \geq 0.
\end{align*}
Using Fourier-Motzkin elimination, we then eliminate the variable $w_{ij}$ to obtain the following system in our original variables $f, x_i,x_j$ and $x_k$.
\begin{align}
&x_i - a_i &\geq 0, \label{in1}\\
&x_j - a_j &\geq 0, \label{in2}\\
&f - a_ja_kx_i - a_ia_kx_j - a_ia_jx_k + 2a_ia_ja_k &\geq 0, \label{in3} \\
&f - a_jb_kx_i - a_ib_kx_j - b_ib_jx_k + a_ia_jb_k + b_ib_jb_k &\geq 0,\label{in4}\\
&- x_j + b_j &\geq 0, \label{in5}\\
&- x_i + b_i &\geq 0, \label{in6}\\
&f - b_ja_kx_i - b_ia_kx_j - a_ia_jx_k + a_ia_ja_k + b_ib_ja_k &\geq 0, \label{in7}\\
&f - b_jb_kx_i - b_ib_kx_j - b_ib_jx_k + 2b_ib_jb_k &\geq 0, \label{in8}\\
&-f + b_jb_kx_i + a_ib_kx_j + a_ia_jx_k - a_ia_jb_k - a_ib_jb_k &\geq 0, \label{in9}\\
&-f + a_jb_kx_i + b_ib_kx_j + a_ia_jx_k - a_ia_jb_k - b_ia_jb_k &\geq 0, \label{in10}\\
&- x_k + b_k &\geq 0, \label{in11}\\
&-f + b_ja_kx_i + a_ia_kx_j + b_ib_jx_k - a_ib_ja_k - b_ib_ja_k &\geq 0, \label{in12}\\
&-f + a_ja_kx_i + b_ia_kx_j + b_ib_jx_k - b_ia_ja_k - b_ib_ja_k &\geq 0, \label{in13}\\
&x_k - a_k &\geq 0, \label{in14} \\
&f - a_ia_jx_k & \geq 0, \label{in15} \\
&-f + b_ib_jx_k & \geq 0. \label{in16}
\end{align}

It is easy to see that the inequalities \ref{in15} and \ref{in16} are redundant:  \ref{in15} is $a_ja_k(\ref{in1}) + a_ia_k(\ref{in2}) + (\ref{in3})$, and \ref{in16} is $b_ja_k(\ref{in6}) + a_ia_k(\ref{in5}) +(\ref{in12})$.

We use the following notation in what follows.  For $i=1,2,3$, \emph{system} $i$
is defined to be the system of inequalities obtained by first grouping the
pair of variables $x_j$ and $x_k$, with $j$ and $k$ different from $i$.
$\Po_i$ is defined to be the solution set of this system.


\subsection{Hull.}
\label{hull}

As we noted earlier, a convex-hull representation for trilinear monomials is known.  From \cite{Meyer04a},
for any labeling that satisfies  \ref{Omega} (or even just the first inequality of \ref{Omega}),
this inequality system which we refer to as $H$ is:
\begin{align}
f - a_2a_3x_1 - a_1a_3x_2 - a_1a_2x_3 + 2a_1a_2a_3 \geq0,& \label{ineq1}\\
f - b_2b_3x_1 - b_1b_3x_2 - b_1b_2x_3 + 2b_1b_2b_3 \geq0,& \label{ineq2}\\
f - a_2b_3x_1 - a_1b_3x_2 - b_1a_2x_3 + a_1a_2b_3 + b_1a_2b_3 \geq0,& \label{ineq3}\\
f - b_2a_3x_1 - b_1a_3x_2 - a_1b_2x_3 + b_1b_2a_3 + a_1b_2a_3 \geq0,& \label{ineq4}\\
f - {\scriptstyle\frac{\eta_1}{b_1-a_1}}x_1 - b_1a_3x_2 - b_1a_2x_3 + \Big({\scriptstyle\frac{\eta_1a_1}{b_1-a_1}} + b_1b_2a_3 + b_1a_2b_3 - a_1b_2b_3\Big) \geq0,&  \label{ineq5}\\
f - {\scriptstyle\frac{\eta_2}{a_1-b_1}}x_1 - a_1b_3x_2 - a_1b_2x_3 + \Big({\scriptstyle\frac{\eta_2b_1}{a_1-b_1}} + a_1a_2b_3 + a_1b_2a_3 - b_1a_2a_3\Big) \geq0,& \label{ineq6}\\
-f + a_2a_3x_1 + b_1a_3x_2 + b_1b_2x_3 - b_1b_2a_3 - b_1a_2a_3 \geq0,& \label{ineq7}\\
-f + b_2a_3x_1 + a_1a_3x_2 + b_1b_2x_3 - b_1b_2a_3 - a_1b_2a_3 \geq0,& \label{ineq8}\\
-f + a_2a_3x_1 + b_1b_3x_2 + b_1a_2x_3 - b_1a_2b_3 - b_1a_2a_3 \geq0,& \label{ineq9}\\
-f + b_2b_3x_1 + a_1a_3x_2 + a_1b_2x_3 - a_1b_2b_3 - a_1b_2a_3 \geq0,& \label{ineq10}\\
-f + a_2b_3x_1 + b_1b_3x_2 + a_1a_2x_3 - b_1a_2b_3 - a_1a_2b_3 \geq0,& \label{ineq11}\\
-f + b_2b_3x_1 + a_1b_3x_2 + a_1a_2x_3 - a_1b_2b_3 - a_1a_2b_3 \geq0,& \label{ineq12}\\
x_1 - a_1 \geq 0,& \label{ineq13}\\
- x_1 + b_1 \geq 0,& \label{ineq14}\\
x_2 - a_2 \geq 0,& \label{ineq15}\\
- x_2 + b_2 \geq 0,& \label{ineq16}\\
x_3 - a_3 \geq 0,& \label{ineq17}\\
- x_3 + b_3 \geq 0,& \label{ineq18}
\end{align}
where $\eta_1 = b_1b_2a_3-a_1b_2b_3-b_1a_2a_3+b_1a_2b_3$ and $\eta_2 = a_1a_2b_3 - b_1a_2a_3 - a_1b_2b_3 + a_1b_2a_3$.\\

We refer to the polytope defined as the feasible set of system $H$ as $\Po_{H}$.   The extreme points of $\Po_{H}$ are the 8 points that correspond to the $2^3=8$ choices of each $x$-variable at its upper or lower bound. We label these 8 points (all of the form $[f=x_1x_2x_3, x_1, x_2, x_3]$) as follows:

{\small\begin{align*}
&v^1:= \left [\begin{array} {c}
b_1a_2a_3\\
b_1\\
a_2\\
a_3
  \end{array} \right] ,\; v^2:= \left [\begin{array} {c}
a_1a_2a_3\\
a_1\\
a_2\\
a_3
  \end{array} \right], \;v^3:= \left [\begin{array} {c}
a_1a_2b_3 \\
a_1\\
a_2\\
b_3
  \end{array} \right], \;v^4:= \left [\begin{array} {cc}
a_1b_2a_3 \\
a_1\\
b_2\\
a_3
  \end{array} \right], \\[5pt] &v^5:= \left [\begin{array} {cc}
a_1b_2b_3 \\
a_1\\
b_2\\
b_3
  \end{array} \right], \;v^6:= \left [\begin{array} {cc}
b_1b_2b_3\\
b_1\\
b_2\\
b_3
  \end{array} \right], \;v^7:=  \left [\begin{array} {cc}
b_1b_2a_3 \\
b_1\\
b_2\\
a_3
  \end{array} \right], \;v^8:= \left [\begin{array} {cc}
b_1a_2b_3 \\
b_1\\
a_2\\
b_3
  \end{array} \right].
  \end{align*}}

Each alternative polyhedral convexification leads to a different system of inequalities (system $i$, $i=1,2,3$)
and therefore a different polytope ($\Po_i$, $i=1,2,3$) in $\R^4$ --- all three contain the convex hull of the solution set of our original trilinear monomial (on the box domain), i.e. $\Po_H$.

To establish if one of these three convexifications is better than another, we need to be able to compare these polytopes in a quantifiable manner. We take the (4-dimensional) volume as our measure, with the idea that a smaller volume corresponds to a tighter convexification. 

For trilinear monomials with domain being a box (in the nonnegative orthant), we derive exact expressions for the  (4-dimensional) volume for the convex hull of the set of solutions and also for each of the three possible double-McCormick convexifications. These volumes are in terms of six parameters (the upper and lower bounds on each of the three variables) and are rather complicated.  By comparing the volume expressions, we are able to draw conclusions regarding the optimal way to perform double McCormick for trilinear monomials.

%% file: alternatives.tex
\section{Alternatives.}\label{sec:alternatives}

In practice, there are many possibilities for handling each product of three
terms encountered in a formulation. A good choice, which may well be different for
different triple products in the same formulation, ultimately depends on trading off
the tightness of a relaxation with the overhead in working with it. For clarity,
in the remainder of this section, we focus on different possible treatments of $f=x_1x_2x_3$.

One possibility is to use the full trilinear hull $\Po_H$. This representation has the
benefit of using no auxiliary variables. Another possibility to use the
\emph{convex-hull representation} (see \cite{Costa}, for example),
writing $f=\sum_{j=1}^8 \lambda_j v^j$, with $\sum_{j=1}^8 \lambda_j =1$,
$\lambda_j\geq 0$,
for $j=1,2,\ldots,8$. This formulation has the drawback of utilizing \emph{eight auxiliary variables}.
But noticing that there are 5 linear equations, we can really reduce to \emph{three auxiliary variables}.
In fact, there is a very structured way to do this, where none of the $\lambda_j$ variables
are employed at all, and rather we introduce  \emph{three auxiliary variables} $w_{12}$,
$w_{13}$ and $w_{23}$, which represent the products $x_1x_2$, $x_1x_3$ and $x_2x_3$, respectively.
A strong advantage of this last approach is when terms $x_1x_2$, $x_1x_3$ and $x_2x_3$
are also in the model under consideration. We wish to emphasize that projecting any of these
convex-hull representations (reduced or not) down to the space of $(f,x_1,x_2,x_3)$
yields again $\Po_H$, and so all of these representations have the same bounding power.

We are advocating the \emph{consideration} of double-McCormick relaxations as an alternative
when warranted.
We have identified the best among the double McCormicks and quantified the
error in using it in preference to $\Po_H$ (and, ipso facto, with any convex-hull
or reduced convex hull representation).
A double-McCormick relaxation involves only \emph{one auxiliary variable} (and 8 inequalities).
This can be particularly attractive
when this particular auxiliary variable already appears in the model under consideration.
Alternatively, especially when this particular auxiliary variable does \emph{not} appear in the formulation,
we can use the formulation with \emph{zero auxiliary variables} (\ref{in1}-\ref{in14}).
Recently (see \cite{SpeakmanLee2016b}), we have computationally validated such an approach in the context of ``box cubic programs''
\begin{equation*}
\label{boxcup}
 \min_{x\in\mathbb{R}^n}\left\{ \sum_{\{i,j,k\}} q_{\tiny ijk}\,   x_i\,   x_j\,  x_k ~:~
x_i \in[a_i,b_i],~ i=1,2,\ldots,n\right\}. 
\end{equation*}
In this type of problem, we can apply (\ref{in1}-\ref{in14}) \emph{independently} for each trinomial, with no
auxiliary variables at all, choosing the best double-McCormick for each trinomial, whenever the associated volume is close
to the volume for $\Po_H$. We have documented that this can happen quite a lot, and so it is a viable approach.
It is important to emphasize that some of the negative experience with double McCormick is
related to choosing the wrong one. Indeed, our mathematical and computational results indicate
that there are many situations where: (i) the worst double McCormick is quite bad compared to the best one,
and (ii) the best one is only slightly worse than $\Po_H$ (and its convex-hull representations).

Besides any prescriptive use of double-McCormick relaxations, our results can simply be seen as
quantifying the bounding advantage given by $\Po_H$ and the various convex-hull representations (reduced or not)
as compared to each of the possible double McCormick relaxations.

%

In some global-optimization software (e.g., \verb;BARON; and \verb;ANTIGONE;)
the complicated inequality description of the trilinear hull is explicitly used.
In other global-optimization software (e.g., \verb;COUENNE; and \verb;SCIP;)
and as a technique at the formulation level,
repeated McCormick is used for the trilinear case. It is by no means
clear that either approach should be followed all of the time
(though this currently seems to be the case), because
of the solution-time tradeoff in using more complicated but stronger
convexifications. This effect can be especially
pronounced in the case of nonlinear optimization where
solutions may not be on the boundary (see \cite{Lee_2007}, for example).
By quantifying the quality of different convexifications, we
offer (i) firm and actionable means for deciding between them at  run time
and, (ii) some explanation for differing behavior of sBB software
under different scenarios.

Finally,
 we note that
the double-McCormick
approach is often applied at the modeling level (see \cite{LejeuneMargot2016}, for example).
In particular, our results are highly relevant to modelers who
simply use global-optimization software, often through a modeling language.
An uniformed modeler can  defeat clever software.
In such a case, it is very
useful for the user to know which double McCormick to use, because
a bad one may negatively affect sBB performance, and all sBB software
that we know will not capture \emph{implicit} triple products in  formulations.


%% file: theorems.tex
\section{Theorems.}\label{sec:theorems}
First we define the following twelve points in $\R^4$, where
$j:= i+1~(\mbox{mod}~3)$ and $k:= i+2~(\mbox{mod}~3)$:

{\footnotesize
\begin{align*} &v_1^9:= \left [\begin{array} {cc}
\theta_1^1 \\
\theta_1^2\\
a_2\\
b_3
  \end{array} \right],  v_1^{10}:= \left [\begin{array} {cc}
\theta_1^3 \\
\theta_1^4\\
b_2\\
a_3
  \end{array} \right],  v_1^{11}:= \left [\begin{array} {cc}
\theta_1^5 \\
\theta_1^6\\
b_2\\
a_3
  \end{array} \right], v_1^{12}:= \left [\begin{array} {cc}
\theta_1^7 \\
\theta_1^8\\
a_2\\
b_3
  \end{array} \right],
  v_2^9:= \left [\begin{array} {cc}
\theta_2^1 \\
b_1\\
\theta_2^2\\
a_3
  \end{array} \right],  v_2^{10}:= \left [\begin{array} {cc}
\theta_2^3 \\
a_1\\
\theta_2^4\\
b_3
  \end{array} \right], \\[5pt] &v_2^{11}:= \left [\begin{array} {cc}
\theta_2^5 \\
a_1\\
\theta_2^6\\
b_3
  \end{array} \right],  v_2^{12}:= \left [\begin{array} {cc}
\theta_2^7 \\
b_1\\
\theta_2^8\\
a_3
  \end{array} \right],
v_3^9:= \left [\begin{array} {cc}
\theta_3^3 \\
b_1\\
a_2\\
\theta_3^4
  \end{array} \right],  v_3^{10}:= \left [\begin{array} {cc}
\theta_3^1 \\
a_1\\
b_2\\
\theta_3^2
  \end{array} \right],  v_3^{11}:= \left [\begin{array} {cc}
\theta_3^7 \\
a_1\\
b_2\\
\theta_3^8
  \end{array} \right],  v_3^{12}:= \left [\begin{array} {cc}
\theta_3^5 \\
b_1\\
a_2\\
\theta_3^6
  \end{array} \right],
  \end{align*}}
where:
\begin{align*}
\theta_i^1 &= a_ia_ja_k + \frac{a_j(b_k-a_k)(b_ib_jb_k-a_ia_ja_k)}{b_jb_k-a_ja_k},\;\; &\theta_i^2 = a_i + \frac{a_j(b_i-a_i)(b_k-a_k)}{b_jb_k-a_ja_k},\\
\theta_i^3 &= a_ia_ja_k + \frac{a_k(b_j-a_j)(b_ib_jb_k-a_ia_ja_k)}{b_jb_k-a_ja_k},\;\; &\theta_i^4 = a_i + \frac{a_k(b_j-a_j)(b_i-a_i)}{b_jb_k-a_ja_k},\\
\theta_i^5 &= \frac{b_ja_k(a_ib_jb_k-a_ia_jb_k-b_ia_ja_k+b_ia_jb_k)}{b_jb_k-a_ja_k},\;\; &\theta_i^6 = a_i + \frac{b_j(b_i-a_i)(b_k-a_k)}{b_jb_k-a_ja_k},\\
\theta_i^7 &= \frac{a_jb_k(b_ib_ja_k-b_ia_ja_k-a_ib_ja_k+a_ib_jb_k)}{b_jb_k-a_ja_k},\;\; &\theta_i^8 = a_i + \frac{b_k(b_j-a_j)(b_i-a_i)}{b_jb_k-a_ja_k}.
\end{align*}

Next, we state our main results.
\begin{theorem}
\label{TheoremH}
\begin{dmath*}
\Vol_{\Po_H} = (b_1-a_1)(b_2-a_2)(b_3-a_3)\times\\
\left(b_1(5b_2b_3 - a_2b_3 - b_2a_3 - 3a_2a_3) + a_1(5a_2a_3 - b_2a_3 - a_2b_3 -3b_2b_3)\right)/24.
\end{dmath*}
\end{theorem}

\begin{theorem}
\label{TheoremB}
The set of extreme points of $\Po_1$ is
$\{v^1,\ldots, v^8\}\cup\{v_1^9,\ldots,v_1^{12}\}$.
Moreover,
\begin{dmath*}
\Vol_{\Po_1} = \Vol_{\Po_H} + (b_1-a_1)(b_2-a_2)^2(b_3-a_3)^2\times\\
\frac{3(b_1b_2a_3 - a_1b_2a_3 + b_1a_2b_3 -a_1a_2b_3) + 2(a_1b_2b_3 - b_1a_2a_3)}{24(b_2b_3-a_2a_3)}.
\end{dmath*}
\end{theorem}

\begin{theorem}
\label{TheoremC}
The set of extreme points of $\Po_2$ is
$\{v^1,\ldots, v^8\}\cup\{v_2^9,\ldots,v_2^{12}\}$.
Moreover,
\begin{dmath*}
\Vol_{\Po_2} = \Vol_{\Po_H} +
\frac{(b_1-a_1)(b_2-a_2)^2(b_3-a_3)^2\left(5(a_1b_1b_3-a_1b_1a_3) + 3(b_1^2a_3 - a_1^2b_3)\right)}{24(b_1b_3-a_1a_3)}.
\end{dmath*}
\end{theorem}

\begin{theorem}
\label{TheoremA}
The set of extreme points of $\Po_3$ is
$\{v^1,\ldots, v^8\}\cup\{v_3^9,\ldots,v_3^{12}\}$.
Moreover,
\begin{dmath*}
\Vol_{\Po_3} = \Vol_{\Po_H}+
\frac{(b_1-a_1)(b_2-a_2)^2(b_3-a_3)^2\left(5(a_1b_1b_2-a_1b_1a_2) + 3(b_1^2a_2 -a_1^2b_2)\right)}{24(b_1b_2-a_1a_2)}.
\end{dmath*}
\end{theorem}

Our proofs in \S\S\ref{Three Continuous Convex Hull}--\ref{linC} all assume that $a_1,a_2,a_3>0$.
Next, we briefly explain why the theorems hold even when any of the $a_i$ are zero.
Taking the convex hull of a compact set is continuous (even 1-Lipschitz) in the Hausdorff metric (see \cite[p. 51]{Schneider1}).
The volume functional is continuous (with respect to the Hausdorff metric) on the set $K^n$ of convex bodies in $\mathbb{R}^n$
(see \cite[Theorem 1.8.20; p. 68]{Schneider2}).
If two sets of $m$ points in $\mathbb{R}^n$ are close as vectors in $\mathbb{R}^{mn}$, then they are also close in the Hausdorff metric.
Therefore, the volume of the convex hull of a set of $m$ points in $\mathbb{R}^n$ is a continuous function of the
coordinates of the points. Also, the coordinates of the extreme points of our polytopes are all continuous functions
(of the six parameters) at $a_i=0$.
Finally, we note that the volume formulae that we derive are continuous functions (of the six parameters)
at $a_i=0$. Therefore, those formulae are also correct when some $a_i=0$.
We do note that we can also modify our constructions to handle these cases where some of the $a_i$ are zero, but our continuity argument
is much shorter.
\\

\begin{corollary}
For all values of the parameters $a_1,b_1,a_2,b_2,a_3,b_3$, meeting the conditions \eqref{Omega}, we have: $\Vol_{\Po_H} \leq \Vol_{\Po_3} \leq \Vol_{\Po_2} \leq \Vol_{\Po_1}.$\\
\end{corollary}

From this we can see that with the variables ordered according to their upper and lower bounds per \eqref{Omega}, the smallest volume will always be obtained by using system 3 (i.e., first grouping variables $x_1$ and $x_2$). In addition, for different values of the upper and lower bounds, we can precisely quantify the difference in volume of the alternative convexifications.

Moreover, by substituting $a_1=a_2=0$ and $b_1=b_2=1$ into the conditions \eqref{Omega}, we can easily see the following corollary relevant to mixed-integer nonlinear optimization.\\

\begin{corollary}
In the special case of two binary variables and one continuous variable, first grouping the two binary variables gives the convexification with the smallest volume.\\
\end{corollary}

\begin{figure}
  \caption{Graph of difference in volume $\left(\frac{3a_3(b_3-a_3)^2}{24b_3}\right)$ vs. parameter values}
  \label{fig1}
  \centering
   \includegraphics[scale=0.6]{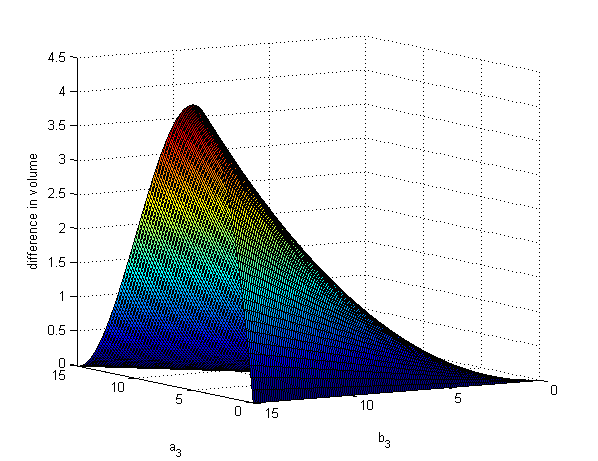}
\end{figure}

In this special case, we only have two parameters $a_3$ and $b_3$ and the volume formulae simplify considerably.  In particular, for this special case, $\Po_3$ is equivalent to $\Po_H$, and $\Po_1$ and $\Po_2$ are equivalent.
We compute the difference in volume between the two distinct choices of convexification and,
 in Figure \ref{fig1}, plot this expression as the parameters vary (satisfying $0\leq a_3 < b_3$).   The following is easy to establish.\\

\begin{corollary}
As $a_3$ and $b_3$ increase, the difference in volumes of $\Po_3$ and $\Po_2$ (or $\Po_1$) becomes arbitrarily large.
Additionally, for a fixed $b_3$, the greatest difference in volume occurs when $a_3=b_3/3$.\\
\end{corollary}

Finally we note that in the special case in which $a_1=a_2=a_3=0$, each convexification reduces to the convex hull, which is a result of \cite{RyooSahi}.
So in this case, \emph{any} double-McCormick convexification has the power of the more-complicated inequality description of the convex hull.
In fact, viewed this way, our results provide a quantified generalization of this
result of \cite{RyooSahi}. We do wish to emphasize that because our results
do not just apply to trilinear monomials on the formulation variables,
but may well involve auxiliary variables, \emph{the case of non-zero lower bounds
is very relevant.}

%% file: theorem_H.tex
\section{Proof of Thm.~\ref{TheoremH}.}
\label{Three Continuous Convex Hull}
\begin{figure}
  \caption{Visual representation of simplex, $\mathcal{S}$, with extreme points $v^1, v^2, v^4, v^5$ and $v^6$, plus the other 3 convex hull extreme points $v^3$, $v^7$ and $v^8$.}
   \label{fig2}
  \centering
      \includegraphics[scale=0.4]{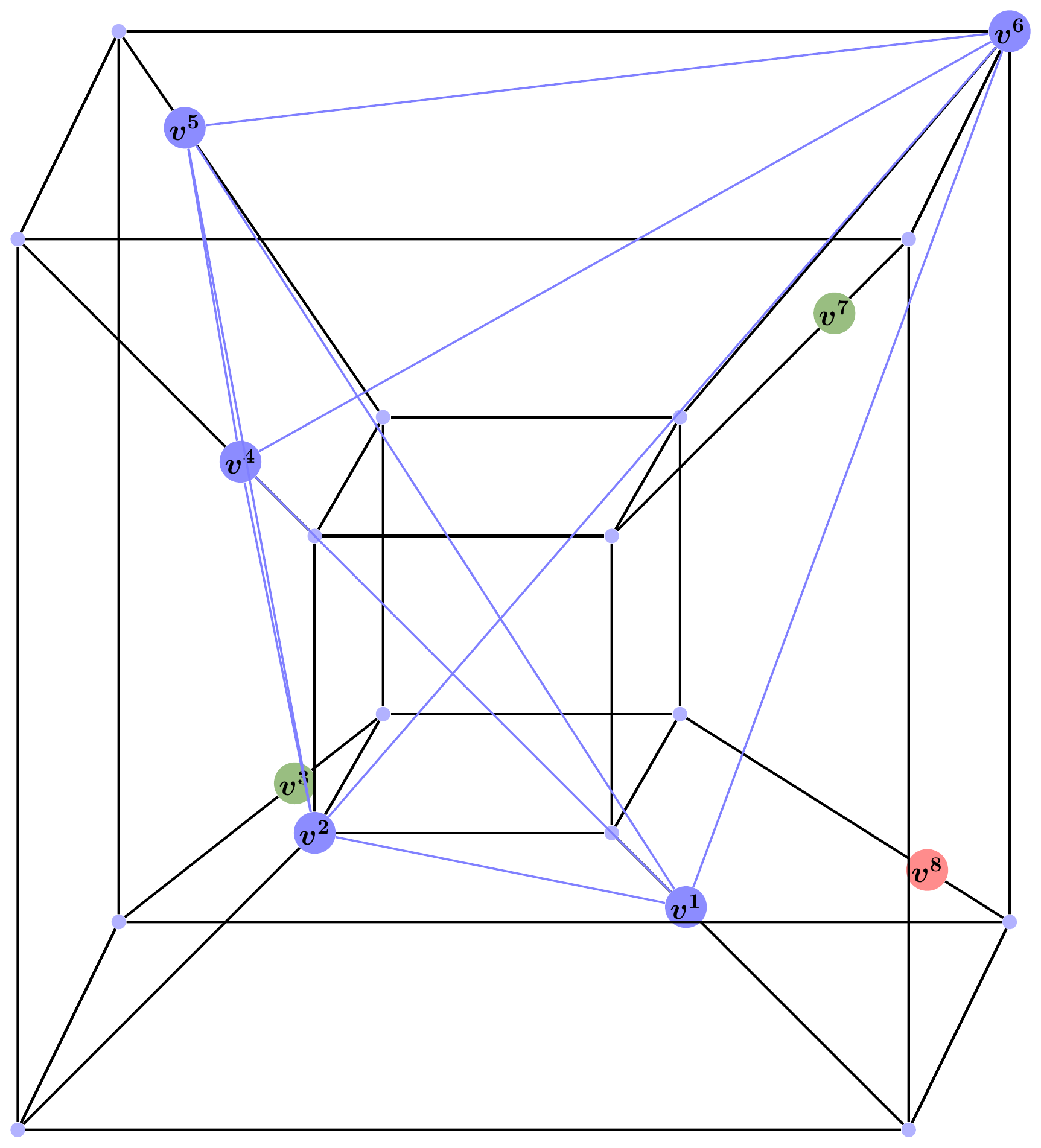}
\end{figure}

We compute the volume of $\Po_H$ by constructing a triangulation.
See Figure \ref{fig2} for a diagram of the 8 extreme points of $\Po_H$.
Note that $v^2$, which has all of the variables at their lower bounds,
is at the bottom of the ``inner cube'', and $v^6$, which has all of the variables at their upper bounds,
is at the top of the ``outer cube''.

We use the fact that
the volume of an $n$-simplex in $\mathbb{R}^n$ with vertices $(z^0,\dots,z^n)$ is:
\[
\big |\det(z^1-z^0 \;\; z^2-z^0 \;\dots\; z^n-z^0)\big |/n!.
\]
The volume of the 4-simplex with extreme points $v^1, v^2, v^4, v^5$ and $v^6$, which we define as $\mathcal{S}\colonequals \text{conv}\{v^1, v^2, v^4, v^5, v^6\}$,
is
\[
{(b_1-a_1)^2(b_2-a_2)(b_3-a_3)(b_2b_3-a_2a_3)}/24.
\]

A 4-simplex has 5 facets, each of which is a 3-simplex and is described by the hyperplane through a choice of 4 extreme points.
To determine the facet-describing inequalities, we compute each hyperplane and then check the final point to obtain the direction of the inequality.  The 5 facets of $\mathcal{S}$ are described as follows:

$F^{1}$ (hyperplane through points $v^1, v^2, v^4, v^6$):
\begin{dmath*}-f + a_2a_3x_1+a_1a_3x_2+\frac{(a_1a_2a_3-a_1b_2a_3-b_1a_2a_3+b_1b_2b_3)}{(b_3-a_3)}x_3-\frac{(a_1a_2a_3b_3-a_1b_2a_3^2-b_1a_2a_3^2+b_1b_2a_3b_3)}{(b_3-a_3)}\geq 0 \end{dmath*}

$F^{2}$ (hyperplane through points $v^1, v^2, v^4, v^5$):
\begin{dmath*}f-a_2a_3x_1-a_1a_3x_2-a_1b_2x_3 +a_1a_2a_3+a_1b_2a_3\geq 0  \end{dmath*}

$F^{3}$ (hyperplane through points $v^1, v^2, v^5, v^6$):
\begin{dmath*}(b_3-a_3)x_2-(b_2-a_2)x_3 + b_2a_3-a_2b_3 \geq 0\end{dmath*}

$F^{4}$ (hyperplane through points $v^1, v^4, v^5, v^6$):
\begin{dmath*}f - b_2b_3x_1-\frac{(a_1b_2a_3-a_1b_2b_3-b_1a_2a_3+b_1b_2b_3)}{(b_2-a_2)}x_2-a_1b_2x_3+\frac{(-a_1a_2b_2b_3+a_1b_2^2a_3-b_1a_2b_2a_3+b_1b_2^2b_3)}{(b_2-a_2)}\geq 0 \end{dmath*}

$F^{5}$ (hyperplane through points $v^2, v^4, v^5, v^6$):
\begin{dmath*}-f+b_2b_3x_1+a_1a_3x_2+a_1b_2x_3-a_1b_2a_3-a_1b_2b_3\geq0\end{dmath*}

If a hyperplane $H$ intersects a polytope $P$ on a facet $F$, then  $H^+$ (resp., $H^-$) denotes the half-space determined by $H$ that contains (does not contain) $P$.  If a point $w$ is not in $H$ but in $H^+$ (resp., $H^-$), then  $w$ is \emph{beneath} (\emph{beyond}) $F$ (see \cite[p. 78]{Grunbaum}).

We now compute the volume of $\text{conv}(\mathcal{S} \cup \{v^8\})$.  To obtain the additional volume of this polytope compared with $\mathcal{S}$ we sum the volume of $\text{conv}(\{v^8\} \cup F)$ for each facet, $F$, of $\mathcal{S}$ such that $v^8$ is beyond that facet.  To do this we first check each of the 5 facets to determine if $v^8$ is beneath or beyond that facet.  To do this we substitute $v^8$ into the relevant inequality, and if the result is negative then $v^8$ lies beyond that facet.


It is easy to check that $v^8$ satisfies $F^1$ and $F^5$ and violates $F^3$.  Using Lemma \ref{lem0}, we also check that $v^8$ satisfies both $F^2$ and $F^4$.  From this we have that $v^8$ is beyond one facet, $F^{3}$.  Therefore we need to calculate the volume of $\text{conv}(F^3 \cup \{v^8\})= \text{conv}\{v^1, v^2, v^5, v^6, v^8\}$, this a a 4-simplex with volume:
\begin{dmath*}{(b_1-a_1)^2(b_2-a_2)(b_3-a_3)(b_2b_3-a_2a_3)}/24.\end{dmath*}

We now have a new polytope which is $\text{conv}\{v^1, v^2, v^4, v^5, v^6, v^8\} = \text{conv}(\mathcal{S} \cup \{v^8\})$.  We refer to this polytope as $\Q$.  The volume of $\Q$ is given by the sum of the volumes of the two simplices we have computed thus far.  The facets of $\Q$ are the facets of the original simplex without $F^{3}$, along with the facets of the 4-simplex: $\text{conv}(F^3 \cup \{v^8\})$ (again not including $F^{3}$ itself).  A facet of $\text{conv}(F^3 \cup \{v^8\})$ is supported by a hyperplane through a choice of 4 of the 5 extreme points (points $v^1, v^2, v^5, v^6$ and $v^8$).  As before, to determine these facet inequalities we compute each hyperplane and then check the final point to obtain the direction of the inequality (note that we exclude the choice $v^1, v^2, v^5, v^6$ because this corresponds to $F^{3}$).  The 4 facets are described below:

$F^{6}$ (plane through points $v^1, v^2, v^5, v^8$):
\begin{dmath*}f-a_2a_3x_1-\frac{(-a_1a_2a_3+a_1b_2b_3+b_1a_2a_3-b_1a_2b_3)}{(b_2-a_2)}x_2-b_1a_2x_3+\frac{(-a_1a_2^2a_3+a_1a_2b_2b_3-b_1a_2^2b_3+b_1a_2b_2a_3)}{(b_2-a_2)}\geq 0\end{dmath*}

$F^{7}$ (plane through points $v^1, v^5, v^6, v^8$):
\begin{dmath*}f-b_2b_3x_1-b_1b_3x_2-b_1a_2x_3+b_1a_2b_3+b_1b_2b_3 \geq 0\end{dmath*}

$F^{8}$ (plane through points $v^2, v^5, v^6, v^8$):
\begin{dmath*}-f+b_2b_3x_1+b_1b_3x_2+\frac{(-a_1a_2a_3+a_1b_2b_3+b_1a_2b_3-b_1b_2b_3)}{(b_3-a_3)}x_3-\frac{(-a_1a_2a_3b_3+a_1b_2b_3^2+b_1a_2b_3^2-b_1b_2a_3b_3)}{(b_3-a_3)}\geq 0\end{dmath*}

$F^{9}$ (plane through points $v^1, v^2, v^6, v^8$):
\begin{dmath*}-f+a_2a_3x_1+b_1b_3x_2+b_1a_2x_3-b_1a_2a_3-b_1a_2b_3\geq 0\end{dmath*}

The facets of $\Q=\text{conv}\{v^1, v^2, v^4, v^5, v^6, v^8\}$ are therefore $F^{1}$, $F^{2}$, $F^{4}$, $F^{5}$, $F^{6}$, $F^{7}$, $F^{8}$ and $F^{9}$.

To obtain the entire volume of $\Po_H$ we need to consider two further extreme points: $v^3$ and $v^7$.  It would be convenient to add these points separately; i.e., compute the additional volume each produces when added to $\Q$, and sum the results.  As the following lemma shows, this will give the correct volume if the intersection of
the line segment between these points and $\Q$ is not empty.

\begin{lemma}
\label{lemlines}
Let $P$ be a convex polytope and let $w_1$ and $w_2$ be points not in $P$.  Let $L(w_1,w_2)$ be the line segment between $w_1$ and $w_2$.
If $L(w_1,w_2) \cap P \not= \emptyset$, then $conv(P,w_1) \cup conv(P,w_2)$ is convex.  Moreover, in this case, $conv(P,w_1,w_2) = conv(P,w_1) \cup conv(P,w_2)$.
\end{lemma}

\proof{Proof.}
First we show that $conv(P,w_1) \cup conv(P,w_2)$ is convex.  If we show that $L(w_1,w_2)$ is completely contained in $conv(P,w_1) \cup conv(P,w_2)$, then we will be done.
Choose $z \in L(w_1,w_2) \cap P$.
Now consider $L(w_1,z)$. Because $z\in P$, this whole line segment must be in $conv(P,w_1)$.  Similarly consider $L(z,w_2)$; this whole line segment must be contained in $conv(P,w_2)$.  Therefore the whole line segment $L(w_1,w_2)$ must be contained in $conv(P,w_1) \cup conv(P,w_2)$ and therefore this set is convex.

Next, we demonstrate that $conv(P,w_1,w_2) = conv(P,w_1) \cup conv(P,w_2)$. First, choose $y \in conv(P,w_1) \cup conv(P,w_2)$; therefore $y \in conv(P,w_1)$ or $y \in conv(P,w_2)$ (or both); in either case it is clear that $y \in conv(P,w_1,w_2)$.  In the other direction, choose $y \in conv(P,w_1,w_2)$; therefore $y$ can be written as a convex combination of the extreme points of $P$ and $w_1$ and $w_2$.  Because $conv(P,w_1) \cup conv(P,w_2)$ is convex, this means $y \in conv(P,w_1) \cup conv(P,w_2)$.  Therefore the sets are equal as required. \halmos
\endproof

We refer to the midpoint of the line between $w_1$ and $w_2$ as $M(w_1,w_2)$.
To show that the intersection of $L(v^3,v^7)$ and $\Q$ is non-empty, consider the midpoint
\begin{equation*} 
\label{v3v7mid} M(v^3,v^7)=\left [\begin{array} {cccc} \frac{a_1a_2b_3 + b_1b_2a_3}{2} &\frac{a_1+b_1}{2}&\frac{a_2+b_2}{2}&\frac{b_3+a_3}{2}\end{array} \right].
\end{equation*}

We show that this point satisfies each of the inequalities of $\Q$ by substituting into each inequality and checking the result.  By showing that each resulting quantity is nonnegative, we  conclude that the midpoint intersects $\Q$. It is easy to see that the midpoint $M(v^3,v^7)$ satisfies $F^1$, $F^5$, $F^8$ and $F^9$.  Using Lemma \ref{lem0}, we also check that $M(v^3,v^7)$ satisfies $F^2$, $F^4$, $F^6$ and $F^7$.  Therefore $\text{conv}(\Q \cup \{v^3\}) \cup \text{conv}(\Q \cup \{v^7\}) = \text{conv}(\Q \cup \{v^3\} \cup \{v^7\}) = \Po_H$.


\subsection[4.1]{Computing the (additional) volume of $\text{conv}\bm{(\Q \cup \{v^3\})}$.}

We now compute the additional volume of $\text{conv}(\Q \cup \{v^3\})$ compared to the volume of $\Q$.  To obtain this, we sum the volumes of $\text{conv}(\{v^3\} \cup F)$ for each facet, $F$, of $\Q$ such that $v^3$ is beyond that facet.  We substitute $v^3$ into each relevant inequality, and if the result is negative then $v^3$ lies beyond that facet.  It is easy to see that $v^3$ satisfies $F^5$, $F^7$ and $F^9$ and violates $F^2$, $F^6$ and $F^8$.  It can then be checked that $v^3$ satisfies $F^1$ using Lemma \ref{lem3} (with $A=b_2, B=a_2, C=(b_1b_3-a_1a_3), D=(2a_1a_3-a_1b_3-b_1a_3)$).   We also check that $v^3$ satisfies $F^4$ using Lemma \ref{lem3} (with $A=b_2, B=a_2, C=(a_1a_3-2a_1b_3+b_1b_3), D=(a_1b_3-b_1a_3)$) and Lemma \ref{lem0}.


From this we know that $v^3$ is beyond $F^{2}$, $F^{6}$ and $F^{8}$; therefore we need to compute the volume of the convex hulls of $v^3$ with each of these facets.

The polytope $\text{conv}(F^2 \cup \{v^3\})= \text{conv}\{v^1, v^2, v^4, v^5, v^3\}$ is a 4-simplex with volume:
\begin{dmath*}{a_1(b_1-a_1)(b_2-a_2)^2(b_3-a_3)^2}/24.\end{dmath*}

The polytope $\text{conv}(F^6 \cup \{v^3\})= \text{conv}\{v^1, v^2, v^5, v^8, v^3\}$ is a 4-simplex with volume:
\begin{dmath*}{a_2(b_1-a_1)^2(b_2-a_2)(b_3-a_3)^2}/24.\end{dmath*}

The polytope $\text{conv}(F^8 \cup \{v^3\})= \text{conv}\{v^2, v^5, v^6, v^8, v^3\}$ is a 4-simplex with volume:
\begin{dmath*}{b_3(b_1-a_1)^2(b_2-a_2)^2(b_3-a_3)}/24.\end{dmath*}

\subsection[4.2]{Computing the (additional) volume of $\text{conv}\bm{(\Q \cup \{v^7\})}$.}

We now compute the additional volume of $\text{conv}(\Q \cup \{v^7\})$ compared to the volume of $\Q$.  To obtain this, we sum the volumes of $\text{conv}(\{v^7\} \cup F)$ for each facet, $F$, of $\Q$ such that $v^7$ is beyond that facet.  We substitute $v^7$ into each relevant inequality, and if the result is negative then $v^7$ lies beyond that facet.  It is easy to see that $v^7$ satisfies $F^2$, $F^5$ and $F^9$ and violates $F^1$, $F^4$ and $F^7$.  It can then be checked that $v^7$ satisfies $F^6$ using Lemma \ref{lem3} (with $A=b_2, B=a_2, C=(b_1a_3-a_1b_3), D=(a_1a_3-2b_1a_3+b_1b_3)$) and Lemma \ref{lem0}.   We  also check that $v^7$ satisfies $F^8$ using Lemma \ref{lem3} (with $A=b_2, B=a_2, C=(2b_1b_3-a_1b_3-b_1a_3), D=(a_1a_3-b_1b_3)$).


From this we know that $v^7$ is beyond $F^{1}$, $F^{4}$ and $F^{7}$, therefore we need to compute the volume of the convex hulls of $v^7$ with each of these facets.

The polytope $\text{conv}(F^1 \cup \{v^7\})= \text{conv}\{v^1, v^2, v^4, v^6, v^7\}$ is a 4-simplex with volume:
\begin{dmath*}{a_3(b_1-a_1)^2(b_2-a_2)^2(b_3-a_3)}/24.\end{dmath*}

The polytope $\text{conv}(F^4 \cup \{v^7\})= \text{conv}\{v^1, v^4, v^5, v^6, v^7\}$ is a 4-simplex with volume:
\begin{dmath*}{b_2(b_1-a_1)^2(b_2-a_2)(b_3-a_3)^2}/24.\end{dmath*}

The polytope $\text{conv}(F^7 \cup \{v^7\})= \text{conv}\{v^1, v^5, v^6, v^8, v^7\}$ is a 4-simplex with volume:
\begin{dmath*}{b_1(b_1-a_1)(b_2-a_2)^2(b_3-a_3)^2}/24.\end{dmath*}

To compute the volume of $\Po_H$, we  sum the volume of the appropriate eight simplices,
and we obtain the volume of $\Po_H$ as stated in Theorem \ref{TheoremH}.
\halmos

%% file: theorem_A.tex
\section{Proof of Thm.~\ref{TheoremA}.}
\label{linA}
We compute the volume of the convex hull of the 12 extreme points which we claim are exactly the extreme points of system 3.  In computing the volume of this polytope, we also prove that these are the correct extreme points and that we have therefore computed the volume of $\Po_3$.

The relevant points are the eight extreme points of $\Po_H$, plus an additional four points.  Because we have already computed the volume of $\Po_H$, to compute the volume of $\Po_3$, we need to compute the additional volume, compared with $\Po_H$, added by these four extra extreme points.  To show that this is indeed the volume of $\Po_3$, we keep track of which facets need to be deleted and added to the system of inequalities as we go.  When this is complete, we have exactly system 3, and therefore we must also have the correct extreme points.

We begin with system $H$ from \S\ref{hull}.  As discussed in \S\ref{Three Continuous Convex Hull}, it would be convenient to add the four points to $\Po_H$ separately; i.e., compute the additional volume each produces when added to $\Po_H$, and sum the results.  To show that we can add two points separately and obtain the correct volume, we show that the intersection of the line segment between these points and $\Po_{H}$ is non-empty (Lemma \ref{lemlines}).

We show that we can add $v_{3}^9$ separately, $v_{3}^{10}$ separately, and then $v_{3}^{11}$ and $v_{3}^{12}$ together by considering the midpoints of the line segments between the relevant points.  We consider $L(v_{3}^9,v_{3}^{10})$, $L(v_{3}^9,v_{3}^{11})$,
$L(v_{3}^9,v_{3}^{12})$, $L(v_{3}^{10},v_{3}^{11})$ and $L(v_{3}^{10},v_{3}^{12})$.
We show that the midpoint of each line segment satisfies each of the inequalities of $\Po_{H}$ by substituting this point into each inequality and checking the result.
See Table \ref{LinALines} for a summary of the resulting substitutions.  The table notes whether nonnegativity of the resulting quantity follows immediately (after factoring), or by use of a technical lemma (after further explanation in the appendix), or after being rewritten in the way referenced in Figure \ref{LinAlinestablefig}.  Because we have shown that each resulting quantity is nonnegative, we conclude that each of the midpoints intersect $\Po_{H}$, and therefore we can add $v_{3}^9$ separately, $v_{3}^{10}$ separately, and then $v_{3}^{11}$ and $v_{3}^{12}$ together.

\begin{table}
\caption{}\label{LinALines}
{\footnotesize  \begin{center}
  \begin{tabular}{ l | c | c |c |c |c }
     Ineq & $M(v_3^{9}, v_3^{10})$ & $M(v_3^{9}, v_3^{11})$ & $M(v_3^{9}, v_3^{12})$ & $M(v_3^{10}, v_3^{11})$ & $M(v_3^{10}, v_3^{12})$\\ \hline
    \ref{ineq1} & immediate & immediate & immediate & immediate & immediate\\
   \ref{ineq2} & immediate & immediate & immediate & immediate & immediate\\
  \ref{ineq3} & immediate & by Lemma \ref{lem0} & immediate & by Lemma \ref{lem0} & by Lemma \ref{lem0}\\
  \ref{ineq4} & immediate & by Lemma \ref{lem0} & by Lemma \ref{lem0} & immediate & by Lemma \ref{lem0}\\
  \ref{ineq5} & immediate & by Lemma \ref{lem0} & immediate & by Lemma \ref{lem0} & by Lemma \ref{lem0}\\
 \ref{ineq6} & immediate & by Lemma \ref{lem0} & by Lemma \ref{lem0} & immediate & by Lemma \ref{lem0}\\
  \ref{ineq7} & immediate & immediate & immediate & immediate & immediate\\
   \ref{ineq8} & immediate & immediate & immediate & immediate & immediate\\
   \ref{ineq9} & see \ref{table:midpoint:1} & see \S\ref{A:midpoint911:ineq9} & immediate & see \ref{table:midpoint:2} & see \ref{table:midpoint:3} \\
    \ref{ineq10} & see \ref{table:midpoint:4}  & see \ref{table:midpoint:5}  & see \ref{table:midpoint:2} & immediate & See \S\ref{A:midpoint1012:ineq10}\\
 \ref{ineq11} & immediate & immediate & immediate & immediate & immediate\\
  \ref{ineq12} & immediate & immediate & immediate & immediate & immediate\\
   \ref{ineq13} & immediate & immediate & immediate & immediate & immediate\\
 \ref{ineq14} & immediate & immediate & immediate & immediate & immediate\\
  \ref{ineq15} & immediate & immediate & immediate & immediate & immediate\\
  \ref{ineq16} & immediate & immediate & immediate & immediate & immediate\\
  \ref{ineq17} & immediate & immediate & immediate & immediate & immediate\\
  \ref{ineq18} & immediate & immediate & immediate & immediate & immediate\\
    \hline
  \end{tabular}
\end{center}}
\end{table}

\begin{figure}
\caption{For Table \ref{LinALines}}
{\footnotesize \label{LinAlinestablefig}
\begin{equation}
\frac{(b_2-a_2)(b_1-a_1)\left (b_1b_3(b_2-a_2)+a_2a_3(b_1-a_1)\right)}{2(b_1b_2-a_1a_2)} \label{table:midpoint:1}
\end{equation}
\begin{equation}
\frac{(b_2b_3-a_2a_3)(b_1-a_1)+(b_1b_3-a_1a_3)(b_2-a_2)}{2}\label{table:midpoint:2}
\end{equation}
\begin{equation}
\frac{(b_2-a_2)(b_1-a_1)\left(b_1(b_2b_3-a_2a_3)+a_2(b_1b_3-a_1a_3)\right)}{2(b_1b_2-a_1a_2)}\label{table:midpoint:3}
\end{equation}
\begin{equation}
\frac{(b_2-a_2)(b_1-a_1)(b_2b_3(b_1-a_1)+a_1a_3(b_2-a_2))}{2(b_1b_2-a_1a_2)}\label{table:midpoint:4}
\end{equation}
\begin{equation}
\frac{(b_2-a_2)(b_1-a_1)\left(b_2(b_1b_3-a_1a_3)+a_1(b_2b_3-a_2a_3)\right)}{2(b_1b_2-a_1a_2)}\label{table:midpoint:5}
\end{equation}}
\end{figure}

\subsection[5.1]{Computing the (additional) volume of $\text{conv}\bm{(\Po_H \cup \{v_3^9\})}$.}

We now compute the additional volume of $\text{conv}(\Po_H \cup \{v_3^9\})$ compared to the volume of $\Po_H$.  To do this we sum the volumes of $\text{conv}(\{v_3^9\} \cup F)$ for each facet, $F$, of $\Po_H$ such that $v_3^9$ is beyond that facet.  We substitute $v_{3}^9$ into each  inequality of system $H$, and we   immediately see that it satisfies every inequality  except \ref{ineq9}.

From this we know that $v_{3}^9$ is beyond only one facet.  The extreme points that lie on this facet are points $v^1, v^2, v^6$ and $v^8$.  The polytope $\text{conv}\{v^1, v^2, v^6, v^8, v_{3}^9\}$ is a 4-simplex with volume:
\begin{dmath*}b_1a_2(b_1-a_1)^2(b_2-a_2)^2(b_3-a_3)^2/\left(24(b_1b_2-a_1a_2)\right ).\end{dmath*}

The facets of $\text{conv}(\Po_H \cup \{v_3^9\})$ are the facets of $\Po_H$ except inequality \ref{ineq9}.  We  see this by computing the four additional facets that come from adding $v_{3}^9$ and noting they are already contained in system $H$:
\begin{itemize}
\item The facet through points $v^1$, $v^2$, $v^6$ and $v_3^9$ is \ref{ineq7}.
\item The facet through points $v^1$, $v^2$, $v^8$ and $v_3^9$ is \ref{ineq15}.
\item The facet through points $v^1$, $v^6$, $v^8$ and $v_3^9$ is \ref{ineq14}.
\item The facet through points $v^2$, $v^6$, $v^8$ and $v_3^9$ is \ref{ineq11}.
\end{itemize}

\subsection[5.2]{Computing the (additional) volume of $\text{conv}\bm{(\Po_H \cup \{v_3^{10}\})}$.}

We now compute the additional volume of $\text{conv}(\Po_H \cup \{v_3^{10}\})$ compared to the volume of $\Po_H$. To do this we sum the volumes of $\text{conv}(\{v_3^{10}\} \cup F)$ for each facet, $F$, of $\Po_H$ such that $v_3^{10}$ is beyond that facet. We substitute $v_{3}^{10}$ into each inequality of system $H$, and we   immediately see that every inequality  is satisfied except \ref{ineq10}.

From this we know that $v_{3}^{10}$ is beyond only one facet.  The extreme points that lie on this facet are points $v^2, v^4, v^5$ and $v^6$.  The polytope $\text{conv}\{v^2, v^4, v^5, v^6, v_{3}^{10}\}$ is a 4-simplex with volume:
\begin{dmath*}a_1b_2(b_1-a_1)^2(b_2-a_2)^2(b_3-a_3)^2/\left(24(b_1b_2-a_1a_2)\right).\end{dmath*}

The facets of $\text{conv}(\Po_H \cup \{v_3^{10}\})$ are the facets of $\Po_H$ except inequality \ref{ineq10}.  We   see this by computing the four additional facets that come from adding $v_{3}^{10}$ and noting that they are already contained in system $H$:
\begin{itemize}
\item The facet through points $v^2$, $v^4$, $v^5$ and $v_3^{10}$ is \ref{ineq13}.
\item The facet through points $v^2$, $v^4$, $v^6$ and $v_3^{10}$ is \ref{ineq8}.
\item The facet through points $v^2$, $v^5$, $v^6$ and $v_3^{10}$ is \ref{ineq12}.
\item The facet through points $v^4$, $v^5$, $v^6$ and $v_3^{10}$ is \ref{ineq16}.
\end{itemize}

\subsection[5.3]{Computing the (additional) volume of $\text{conv}\bm{(\Po_H \cup \{v_3^{11}\} \cup \{v_3^{12}\})}$.}

We now compute the additional volume of $\text{conv}(\Po_H \cup \{v_3^{11}\} \cup \{v_3^{12}\})$ compared to the volume of $\Po_H$.  Because $L(v_{3}^{11},v_{3}^{12})$ lies entirely outside of $\Po_3$, we need to add them sequentially.

We first compute the additional volume of $\text{conv}(\Po_H \cup \{v_3^{11}\})$ compared to the volume of $\Po_H$.  As we have done previously, we sum the volumes of $\text{conv}(\{v_3^{11}\} \cup F)$ for each facet, $F$, of $\Po_H$ such that $v_3^{11}$ is beyond that facet.  We substitute $v_{3}^{11}$ into each relevant inequality, and if the result is negative then $v_{3}^{11}$ lies beyond that facet.  It is immediate that $v_{3}^{11}$ violates inequalities \ref{ineq3}--\ref{ineq6} and satisfies inequalities \ref{ineq1}-\ref{ineq2}, \ref{ineq7}-\ref{ineq8} and \ref{ineq10}-\ref{ineq18}.  To see that inequality \ref{ineq9} is also satisfied see \S\ref{A:point11:ineq9}.

Therefore we have that $v_{3}^{11}$ is beyond four facets, and we need to compute the volume of the convex hulls of $v_{3}^{11}$ with each of these facets.

The extreme points that lie on the first facet are points $v^1, v^3, v^5$ and $v^8$.  The polytope $\text{conv}\{v^1, v^3, v^5, v^8, v_{3}^{11}\}$ is a 4-simplex with volume:
\begin{dmath*}a_1b_1(b_1-a_1)(b_2-a_2)^3(b_3-a_3)^2/\left(24(b_1b_2-a_1a_2)\right).\end{dmath*}

The extreme points that lie on the second facet are points $v^1, v^4, v^5$ and $v^7$.  The polytope $\text{conv}\{v^1, v^4, v^5, v^7, v_{3}^{11}\}$ is a 4-simplex with volume:
\begin{dmath*}a_1b_2(b_1-a_1)^2(b_2-a_2)^2(b_3-a_3)^2/\left(24(b_1b_2-a_1a_2)\right).\end{dmath*}

The extreme points that lie on the third facet are points $v^1, v^5, v^7$ and $v^8$.  The polytope $\text{conv}\{v^1, v^5, v^7, v^8, v_{3}^{11}\}$ is a 4-simplex with volume:
\begin{dmath*}a_1b_1(b_1-a_1)(b_2-a_2)^3(b_3-a_3)^2/\left(24(b_1b_2-a_1a_2)\right).\end{dmath*}

The extreme points that lie on the fourth facet are points $v^1, v^3, v^4$ and $v^5$.  The polytope $\text{conv}\{v^1, v^3, v^4, v^5, v_{3}^{11}\}$ is a 4-simplex with volume:
\begin{dmath*}a_1b_2(b_1-a_1)^2(b_2-a_2)^2(b_3-a_3)^2/\left(24(b_1b_2-a_1a_2)\right).\end{dmath*}

We now have a new polytope which is $\text{conv}(\Po_H \cup \{v_3^{11}\})$.  We refer to this polytope as $\T_3$, and we compute the facets of $\T_3$.

We begin with the facets of $\Po_H$ and delete the four facets that $v_{3}^{11}$ violated (\ref{ineq3}-\ref{ineq6}).  Let us call this system $\T_3^-$.  Now consider the four simplices we dealt with when computing the additional volume produced with $v_{3}^{11}$.  Each of these simplices has 5 facets; one of which corresponds to a deleted facet of $\Po_H$.\\

The remaining 4 facets of the first simplex are described by the planes through the following sets of points:
$\{v^1, v^3, v^5, v_{3}^{11}\}$, $\{v^1, v^3, v^8, v_{3}^{11}\}$, $\{v^1, v^5, v^8, v_{3}^{11}\}$ \\and $\{v^3, v^5, v^8, v_{3}^{11}\}$.

The remaining 4 facets of the second simplex are described by the planes through the following sets of points:
$\{v^1, v^4, v^5, v_{3}^{11}\}$, $\{v^1, v^4, v^7, v_{3}^{11}\}$, $\{v^1, v^5, v^7, v_{3}^{11}\}$ \\and $\{v^4, v^5, v^7, v_{3}^{11}\}$.

The remaining 4 facets of the third simplex are described by the planes through the following sets of points:
$\{v^1, v^5, v^7, v_{3}^{11}\}$, $\{v^1, v^5, v^8, v_{3}^{11}\}$, $\{v^1, v^7, v^8, v_{3}^{11}\}$ \\and $\{v^5, v^7, v^8, v_{3}^{11}\}$.

The remaining 4 facets of the fourth simplex are described by the planes through the following sets of points:
$\{v^1, v^3, v^4, v_{3}^{11}\}$, $\{v^1, v^3, v^5, v_{3}^{11}\}$, $\{v^1, v^4, v^5, v_{3}^{11}\}$ \\and $\{v^3, v^4, v^5, v_{3}^{11}\}$.\\

Consider these sixteen facets and exclude the facets that are shared by more than one simplex.  This leaves eight facets.

We can compute these eight facets to obtain the following:
\begin{itemize}
\item The facet through points $v^1$, $v^3$, $v^8$ and $v_3^{11}$ is
\begin{dmath}
\label{ineq19}
\frac{1}{b_1b_2-a_1a_2}\left (-a_1^2a_2^2b_3+a_1^2a_2b_3x_2-a_1b_1a_2^2a_3+a_1b_1a_2^2x_3\\
{+a_1a_2^2b_3x_1+a_1b_1a_2b_2a_3+a_1b_1a_2a_3x_2-a_1b_1a_2b_3x_2-a_1b_1b_2a_3x_2}\\
+b_1^2a_2b_2b_3-b_1^2a_2b_2x_3-b_1a_2b_2b_3x_1-a_1a_2f+b_1b_2f \right )\geq0 .
\end{dmath}
\item The facet through points $v^3$, $v^5$, $v^8$ and $v_3^{11}$ is
\begin{dmath}
\label{ineq20}
f-a_2b_3x_1-a_1b_3x_2-b_1b_2x_3+a_1a_2b_3+b_1b_2b_3 \geq 0.
 \end{dmath}
\item The facet through points $v^1$, $v^4$, $v^7$ and $v_3^{11}$ is
\begin{dmath}
\label{ineq21}
f-b_2a_3x_1-b_1a_3x_2-a_1a_2x_3+a_1a_2a_3+b_1b_2a_3\geq 0.
\end{dmath}
\item The facet through points $v^4$, $v^5$, $v^7$ and $v_3^{11}$ is \ref{ineq16}.
\item The facet through points $v^1$, $v^7$, $v^8$ and $v_3^{11}$ is
\begin{dmath}
\label{ineq22}
\frac{1}{b_1b_2-a_1a_2}\left (-a_1b_1a_2^2a_3+a_1b_1a_2^2x_3+a_1b_1a_2a_3x_2-a_1b_1a_2b_2b_3+a_1a_2b_2b_3x_1+b_1a_2b_2a_3x_1+b_1^2a_2b_2b_3-b_1^2a_2b_2x_3-b_1a_2b_2b_3x_1+b_1^2b_2^2a_3-b_1^2b_2a_3x_2-b_1b_2^2a_3x_1-a_1a_2f+b_1b_2f \right )\geq 0.
\end{dmath}
\item The facet through points $v^5$, $v^7$, $v^8$ and $v_3^{11}$ is \ref{ineq2}.
\item The facet through points $v^1$, $v^3$, $v^4$ and $v_3^{11}$ is \ref{ineq1}.
\item The facet through points $v^3$, $v^4$, $v^5$ and $v_3^{11}$ is \ref{ineq13}.
\end{itemize}

There are four inequalities that are not already contained in system $\T_3^-$, we add these and in doing so obtain the system of inequalities that describes $\T_3 = \text{conv}(\Po_H \cup \{v_{3}^{11}\})$.

We now compute the additional volume of $\text{conv}(\T_3 \cup \{v_3^{12}\})$ compared to the volume of $\T_3$.  As we have done previously, we sum the volumes of $\text{conv}(\{v_3^{12}\} \cup F)$ for each facet, $F$, of $\T_3$ such that $v_3^{12}$ is beyond that facet.  We substitute $v_{3}^{12}$ into each relevant inequality (i.e., the system of inequalities that describes $\T_3$) and if the result is negative then $v_{3}^{12}$ lies beyond that facet.  It is immediately clear that $v_3^{12}$ satisfies inequalities \ref{ineq1}-\ref{ineq2}, \ref{ineq7}-\ref{ineq9},  \ref{ineq11}-\ref{ineq18} and \ref{ineq20}-\ref{ineq21}.  We can also see immediately that  $v_3^{12}$ violates inequalities \ref{ineq19} and \ref{ineq22}.  To see that inequality \ref{ineq10} is also satisfied see \S\ref{A:point11:ineq9}.

Therefore we see that $v_{3}^{12}$ is beyond two facets and we need to compute the volume of the convex hull of $v_{3}^{12}$ with each of these facets.

The extreme points that lie on the first facet are points $v^1, v^3, v^8$ and $v_{3}^{11}$.  The polytope $\text{conv}\{v^1, v^3, v^8, v_{3}^{11}, v_3^{12}\}$ is a 4-simplex with volume:
\begin{dmath*} b_1a_2(b_1-a_1)^2(b_2-a_2)^2(b_3-a_3)^2/\left(24(b_1b_2-a_1a_2)\right).\end{dmath*}

The extreme points that lie on the second facet are points $v^1, v^7, v^8$ and $v_{3}^{11}$.  The polytope $\text{conv}\{v^1, v^7, v^8, v_{3}^{11}, v_3^{12}\}$ is a 4-simplex with volume:
\begin{dmath*} b_1a_2(b_1-a_1)^2(b_2-a_2)^2(b_3-a_3)^2/\left(24(b_1b_2-a_1a_2)\right).\end{dmath*}

We now compute the additional facets; we take the four facets from adding each simplex and delete the facet that repeats.  This leaves us with the following six facets to compute:
\begin{itemize}
\item The facet through points $v^1$, $v^7$, $v^8$ and $v_3^{12}$ is \ref{ineq14}.
\item The facet through points $v^1$, $v^7$, $v_3^{11}$ and $v_3^{12}$ is \ref{ineq21}.
\item The facet through points $v^7$, $v^8$, $v_3^{11}$ and $v_3^{12}$ is \ref{ineq2}.
\item The facet through points $v^1$, $v^3$, $v^8$ and $v_3^{12}$ is \ref{ineq15}.
\item The facet through points $v^1$, $v^3$, $v_3^{11}$ and $v_3^{12}$ is \ref{ineq1}.
\item The facet through points $v^3$, $v^8$, $v_3^{11}$ and $v_3^{12}$ is \ref{ineq20}.
\end{itemize}

By adding and deleting the appropriate facets to and from system $H$, we see that we arrive at system 3.

Therefore, to compute the volume of $\Po_3$, we sum the volume of $\Po_H$ with that of the appropriate eight simplices, and we obtain our result.
\halmos

%% file: theorem_B.tex
\section{Proof of Thm.~\ref{TheoremB}.}
\label{LinB}
As with Theorem \ref{TheoremA},  we compute the volume of the convex hull of the 12 extreme points which we claim are exactly the extreme points of system 1.  In computing the volume of this polytope, we  also prove that these are the correct extreme points and that we have therefore computed the volume of $\Po_{1}$.

The relevant points are the eight extreme points of $\Po_H$, plus an additional four points.  Because we have already computed the volume of $\Po_H$, to compute the volume of $\Po_{1}$ we need to compute the additional volume, compared with $\Po_H$, added by these four extra extreme points.  To show that this is indeed the volume of $\Po_1$, we keep track of which facets need to be deleted and added to the system of inequalities as we go.  When this is complete, we  have exactly system 1 and therefore we must also have the correct extreme points.

We begin with system $H$ which can be found in \S\ref{hull}, and we use the same principles as we used in the previous proof to compute the volume of $\Po_{3}$.

First we argue that we can add $v_1^9$ to $\Po_H$ separately, $v_1^{10}$ to $\Po_H$ separately and then $v_1^{11}$ and $v_1^{12}$ together.  To do this, we  show that the midpoint of the line segment between $v_1^9$ and all other additional points ($v_1^{10}$, $v_1^{11}$ and $v_1^{12}$) intersects $\Po_H$.  We also show this is true for $v_1^{10}$.

As in the previous proof we refer to the midpoint of the line between $v_i^j$ and $v_i^k$ as $M(v_i^j,v_i^k)$ and we show that the midpoint of each line satisfies each of the inequalities of $\Po_{H}$ by substituting this point into each inequality and checking the result.  See Table \ref{LinBLines} for a summary of the resulting substitutions.  The table notes whether nonnegativity of the resulting quantity follows immediately (after factoring), or by using a technical lemma, after further explanation in the appendix or after being rewritten in the way referenced in Figure \ref{LinBlinestablefig}. Because we have shown that each resulting quantity is nonnegative, we know that the midpoint intersects $\Po_{H}$, and therefore we can add $v_{1}^9$ separately, $v_{1}^{10}$ separately, and then $v_{1}^{11}$ and $v_{1}^{12}$ together.

\begin{table}
\caption{}
\label{LinBLines}
{\footnotesize  \begin{center}
  \begin{tabular}{ l | c | c |c |c |c }
     Ineq & $M(v_1^{9}, v_1^{10})$ & $M(v_1^{9}, v_1^{11})$ & $M(v_1^{9}, v_1^{12})$ & $M(v_1^{10}, v_1^{11})$ & $M(v_1^{10}, v_1^{12})$\\ \hline
    \ref{ineq1} & immediate & immediate & immediate & immediate & immediate\\
   \ref{ineq2} & immediate & immediate & immediate & immediate & immediate\\
  \ref{ineq3} & see \S\ref{B:midpoint910:ineq3} & see \ref{table2:midpoint:1} and & immediate & see \ref{table2:midpoint:2} and& see \S\ref{B:midpoint1012:ineq3}\\
                      &               &   Lemma \ref{lem0} &                       &  Lemma \ref{lem0}  &        \\
   \ref{ineq4} &see \S\ref{B:midpoint910:ineq4} & see \S\ref{B:midpoint911:ineq4} & see \ref{table2:midpoint:2} and & immediate & see \ref{table2:midpoint:4} and \\
                       &              &               &        Lemma \ref{lem0}              &                  &   Lemma \ref{lem0}                   \\
    \ref{ineq5} & see \S\ref{B:midpoint910:ineq5} & see \S\ref{B:midpoint911:ineq5} &  by Lemma \ref{lem0} & by Lemma \ref{lem0} & see \S\ref{B:midpoint1012:ineq5}\\
 \ref{ineq6} & see \S\ref{B:midpoint910:ineq6} & see \S\ref{B:midpoint911:ineq6} & by Lemma \ref{lem0} &  by Lemma \ref{lem0} & see \S\ref{B:midpoint1012:ineq6}\\
  \ref{ineq7} & immediate & immediate & immediate & immediate & immediate\\
   \ref{ineq8} & see \ref{table2:midpoint:5} & see \ref{table2:midpoint:6}& see \ref{table2:midpoint:7} & immediate & see \S\ref{B:midpoint1012:ineq8}\\
   \ref{ineq9} &immediate & immediate & immediate & immediate & immediate \\
  \ref{ineq10} & immediate  & immediate  & immediate & immediate & immediate\\
 \ref{ineq11} & see \ref{table2:midpoint:8} & See \S\ref{B:midpoint911:ineq11} & immediate & see \ref{table2:midpoint:7} & see \ref{table2:midpoint:10}\\
  \ref{ineq12} & immediate & immediate & immediate & immediate & immediate\\
  \ref{ineq13} & immediate & immediate & immediate & immediate & immediate\\
 \ref{ineq14} & immediate & immediate & immediate & immediate & immediate\\
 \ref{ineq15} & immediate & immediate & immediate & immediate & immediate\\
  \ref{ineq16} & immediate & immediate & immediate & immediate & immediate\\
 \ref{ineq17} & immediate & immediate & immediate & immediate & immediate\\
  \ref{ineq18} & immediate & immediate & immediate & immediate & immediate\\
    \hline
  \end{tabular}
\end{center}}
\end{table}

\begin{figure}
\caption{For Table \ref{LinBLines}}
 \label{LinBlinestablefig}
{\footnotesize
\begin{equation}
\frac{(b_3-a_3)(b_2-a_2)\left (b_3(b_1a_2-a_1b_2)+a_2(b_1a_3-a_1b_3)\right)}{2(b_2b_3-a_2a_3)} \label{table2:midpoint:1}
\end{equation}
\begin{equation}
\frac{(b_1a_3-a_1b_3)(b_2-a_2)+(b_1a_2-a_1b_2)(b_3-a_3)}{2}\label{table2:midpoint:2}
\end{equation}
\begin{equation}
\frac{(b_3-a_3)(b_2-a_2)\left (b_2(b_1a_3-a_1b_3)+a_3(b_1a_2-a_1b_2)\right)}{2(b_2b_3-a_2a_3)} \label{table2:midpoint:4}
\end{equation}
\begin{equation}
\frac{(b_3-a_3)(b_2-a_2)\left (b_1b_2(b_3-a_3)+a_1a_3(b_2-a_2)\right)}{2(b_2b_3-a_2a_3)}\label{table2:midpoint:5}
\end{equation}
\begin{equation}
\frac{(b_3-a_3)(b_2-a_2)\left (b_2(b_1b_3-a_1a_3)+a_3(b_1b_2-a_1a_2)\right)}{2(b_2b_3-a_2a_3)}\label{table2:midpoint:6}
\end{equation}
\begin{equation}
\frac{(b_1b_3-a_1a_3)(b_2-a_2)+(b_1b_2-a_1a_2)(b_3-a_3)}{2}\label{table2:midpoint:7}
\end{equation}
\begin{equation}
\frac{(b_3-a_3)(b_2-a_2)\left (b_1b_3(b_2-a_2)+a_1a_2(b_3-a_3)\right)}{2(b_2b_3-a_2a_3)}\label{table2:midpoint:8}
\end{equation}
\begin{equation}
\frac{(b_3-a_3)(b_2-a_2)\left (b_3(b_1b_2-a_1a_2)+a_2(b_1b_3-a_1a_3)\right)}{2(b_2b_3-a_2a_3)}\label{table2:midpoint:10}
\end{equation}}
\end{figure}

\subsection[6.1]{Computing the (additional) volume of $\text{conv}\bm{(\Po_H \cup \{v_1^9\})}$.}
We now compute the additional volume of $\text{conv}(\Po_H \cup \{v_1^9\})$ compared to the volume of $\Po_H$.   To do this, we sum the volumes of $\text{conv}(\{v_1^9\} \cup F)$ for each facet, $F$, of $\Po_H$ such that $v_1^9$ is beyond that facet.   We substitute $v_1^9$ into the 18 relevant inequalities (\ref{ineq1}-\ref{ineq18}) and  immediately see that it satisfies \ref{ineq1}-\ref{ineq3}, \ref{ineq7}-\ref{ineq10} and \ref{ineq12}-\ref{ineq18}.  It is also immediate to see that inequality \ref{ineq11} is violated.  To show that the remaining three inequalities are satisfied (\ref{ineq4}, \ref{ineq5} and \ref{ineq6}) see \S\ref{B:point9:ineq4}, \S\ref{B:point9:ineq5} and \S\ref{B:point9:ineq6}.

From this we know that $v_1^{9}$ is beyond only one facet.  The extreme points that lie on this facet are points $v^2, v^3, v^6$ and $v^8$.  The polytope $\text{conv}\{v^2, v^3, v^6, v^8, v_{1}^9\}$ is a 4-simplex with volume:
\begin{dmath*}a_2b_3(b_1-a_1)^2(b_2-a_2)^2(b_3-a_3)^2/\left(24(b_2b_3-a_2a_3)\right).\end{dmath*}

The facets of $\text{conv}(\Po_H \cup \{v_1^9\})$ are the facets of $\Po_H$ except inequality \ref{ineq11}.  We   see this by computing the four additional facets that come from adding $v_{1}^9$ and noting that they are already contained in system $H$:
\begin{itemize}
\item The facet through points $v^2$, $v^3$, $v^6$ and $v_1^9$ is \ref{ineq12}.
\item The facet through points $v^2$, $v^3$, $v^8$ and $v_1^9$ is \ref{ineq15}.
\item The facet through points $v^2$, $v^6$, $v^8$ and $v_1^9$ is \ref{ineq9}.
\item The facet through points $v^3$, $v^6$, $v^8$ and $v_1^9$ is \ref{ineq18}.
\end{itemize}

\subsection[6.2]{Computing the (additional) volume of $\text{conv}\bm{(\Po_H \cup \{v_1^{10}\})}$.}
We now compute the additional volume of $\text{conv}(\Po_H \cup \{v_1^{10}\})$ compared to the volume of $\Po_H$.   To do this, we sum the volumes of $\text{conv}(\{v_1^{10}\} \cup F)$ for each facet, $F$, of $\Po_H$ such that $v_1^{10}$ is beyond that facet.  We substitute $v_1^{10}$ into the 18 relevant inequalities and   immediately see that it satisfies \ref{ineq1}, \ref{ineq2}, \ref{ineq4}, \ref{ineq7} and \ref{ineq9}-\ref{ineq18}.  It is also immediate to see that inequality \ref{ineq8} is violated.  To show that the remaining three inequalities are satisfied (\ref{ineq3}, \ref{ineq5} and \ref{ineq6}) see \S\ref{B:point9:ineq4}, \S\ref{B:point10:ineq5} and \S\ref{B:point10:ineq6}.

From this we know that $v_1^{10}$ is beyond only one facet.  The extreme points that lie on this facet are points $v^2, v^4, v^6$ and $v^7$.  The polytope $\text{conv}\{v^2, v^4, v^6, v^7, v_{3}^{10}\}$ is a 4-simplex with volume:
\begin{dmath*}b_2a_3(b_1-a_1)^2(b_2-a_2)^2(b_3-a_3)^2/\left(24(b_2b_3-a_2a_3)\right).\end{dmath*}

The facets of $\text{conv}(\Po_H \cup \{v_1^{10}\})$ are the facets of $\Po_H$ except inequality \ref{ineq8}.  We   see this by computing the four additional facets that come from adding $v_{1}^{10}$ and noting that they are already contained in system $H$:
\begin{itemize}
\item The facet through points $v^2$, $v^4$, $v^6$ and $v_1^{10}$ is \ref{ineq10}.
\item The facet through points $v^2$, $v^4$, $v^7$ and $v_1^{10}$ is \ref{ineq17}.
\item The facet through points $v^2$, $v^6$, $v^7$ and $v_1^{10}$ is \ref{ineq7}.
\item The facet through points $v^4$, $v^6$, $v^7$ and $v_1^{10}$ is \ref{ineq16}.
\end{itemize}

\subsection[6.3]{Computing the (additional) volume of $\text{conv}\bm{(\Po_H \cup \{v_1^{11}\} \cup \{v_1^{12}\})}$.}

We now compute the additional volume of $\text{conv}(\Po_H \cup \{v_1^{11}\} \cup \{v_1^{12}\})$ compared to the volume of $\Po_H$.  Because $L(v_1^{11},v_1^{12})$ lies entirely outside of $\Po_{1}$, we need to add them sequentially.

We first compute the additional volume of $\text{conv}(\Po_H \cup \{v_1^{11}\})$ compared to the volume of $\Po_H$.  As we have done previously, we sum the volumes  of $\text{conv}(\{v_1^{11}\} \cup F)$ for each facet, $F$, of $\Po_H$ such that $v_1^{11}$ is beyond that facet.  We substitute $v_1^{11}$ into each relevant inequality and if the result is negative then $v_1^{11}$ lies beyond that facet.  It is immediate that $v_1^{11}$ satisfies inequalities \ref{ineq1}, \ref{ineq2}, \ref{ineq7}-\ref{ineq10} and \ref{ineq12}-\ref{ineq18}.  In \S\ref{B:point11:ineq11} we show that \ref{ineq11} is also satisfied.  It is also immediate that $v_1^{11}$ violates the three facets described by \ref{ineq4}-\ref{ineq6}.  We compute the volume of the convex hulls of $v_1^{11}$ with each of these facets.

The extreme points that lie on the first facet are points $v^1, v^4, v^5$ and $v^7$.  The polytope $\text{conv}\{v^1, v^4, v^5, v^7, v_{1}^{11}\}$ is a 4-simplex with volume:
\begin{dmath*} b_2a_3(b_1-a_1)^2(b_2-a_2)^2(b_3-a_3)^2/\left(24(b_2b_3-a_2a_3)\right).\end{dmath*}

The extreme points that lie on the second facet are points $v^1, v^5, v^7$ and $v^8$.  The polytope $\text{conv}\{v^1, v^5, v^7, v^8, v_{1}^{11}\}$ is a 4-simplex with volume:
\begin{dmath*}b_1a_3(b_1-a_1)(b_2-a_2)^3(b_3-a_3)^2/\left(24(b_2b_3-a_2a_3)\right).\end{dmath*}

The extreme points that lie on the third facet are points $v^1, v^3, v^4$ and $v^5$.  The polytope $\text{conv}\{v^1, v^3, v^4, v^5, v_{1}^{11}\}$ is a 4-simplex with volume:
\begin{dmath*}a_1b_2(b_1-a_1)(b_2-a_2)^2(b_3-a_3)^3/\left(24(b_2b_3-a_2a_3)\right).\end{dmath*}

Unlike in system 3, we see immediately that there exists a fourth facet (described by \ref{ineq3}) which, under certain circumstances, $v_1^{11}$ is beyond.  In particular, this is true if and only if $a_1b_2b_3-b_1a_2a_3 > 0$.  Therefore we continue with two cases.

\subsubsection[s1]{Case 1: $\bm{a_1b_2b_3-b_1a_2a_3 > 0}$}

In this case there exists a fourth facet (described by \ref{ineq3}) such that $v_1^{11}$ is beyond this facet.  The extreme points that lie on this fourth facet are points $v^1, v^3, v^5$ and $v^8$.  The polytope $\text{conv}\{v^1, v^3, v^5, v^8, v_{1}^{11}\}$ is a 4-simplex with volume:
\begin{dmath*}(a_1b_2b_3-b_1a_2a_3)(b_1-a_1)(b_2-a_2)^2(b_3-a_3)^2/\left(24(b_2b_3-a_2a_3)\right).\end{dmath*}

We now have a new polytope which is $\text{conv}(\Po_H \cup \{v_3^{11}\})$ (in case 1).  We refer to this polytope as $\T_1^1$, and compute the facets of $\T_1^1$.

We begin with the facets of $\Po_H$ and delete the four facets that $v_1^{11}$ lies beyond (\ref{ineq3}-\ref{ineq6}).  Let us call this system $\T_1^{1-}$.  Now consider the four simplices we dealt with when computing the additional volume produced with $v_1^{11}$.  Each of these simplices has 5 facets; one of which corresponds to a deleted facet of $\Po_H$.\\

The remaining 4 facets of the first simplex are described by the planes through the following sets of points:
$\{v^1, v^4, v^5, v_1^{11}\}$, $\{v^1, v^4, v^7, v_1^{11}\}$, $\{v^1, v^5, v^7, v_1^{11}\}$ \\and $\{v^4, v^5, v^7, v_1^{11}\}$.

The remaining 4 facets of the second simplex are described by the planes through the following sets of points:
$\{v^1, v^5, v^7, v_1^{11}\}$, $\{v^1, v^5, v^8, v_1^{11}\}$, $\{v^1, v^7, v^8, v_1^{11}\}$ \\and $\{v^5, v^7, v^8, v_1^{11}\}$.

The remaining 4 facets of the third simplex are described by the planes through the following sets of points:
$\{v^1, v^3, v^4, v_1^{11}\}$, $\{v^1, v^3, v^5, v_1^{11}\}$, $\{v^1, v^4, v^5, v_1^{11}\}$ \\and $\{v^3, v^4, v^5, v_1^{11}\}$.

The remaining 4 facets of the fourth simplex are described by the planes through the following sets of points:
$\{v^1, v^3, v^5, v_1^{11}\}$, $\{v^1, v^3, v^8, v_1^{11}\}$, $\{v^1, v^5, v^8, v_1^{11}\}$ \\and $\{v^3, v^5, v^8, v_1^{11}\}$.\\

Consider these sixteen facets and exclude the facets that are shared by more than one simplex.  This leaves eight facets.

We  compute these eight facets to obtain the following:
\begin{itemize}
\item The facet through points $v^1$, $v^3$, $v^8$ and $v_1^{11}$ is
\begin{dmath}
\label{ineq19b}
\frac{1}{b_2b_3-a_2a_3}\left (
{-a_1a_2^2a_3b_3+a_1a_2b_2a_3b_3+a_1a_2a_3b_3x_2-a_1b_2a_3b_3x_2}
-b_1a_2^2a_3^2+b_1a_2^2a_3x_3+a_2^2a_3b_3x_1+b_1a_2a_3^2x_2-b_1a_2a_3b_3x_2+b_1a_2b_2b_3^2-b_1a_2b_2b_3x_3-a_2b_2b_3^2x_1-a_2a_3f+b_2b_3f \right )\geq0 .
\end{dmath}
\item The facet through points $v^3$, $v^5$, $v^8$ and $v_1^{11}$ is
\begin{dmath}
\label{ineq20b}
\frac{1}{b_2b_3-a_2a_3}\left (-a_1a_2^2a_3b_3+a_1a_2a_3b_3x_2+a_1a_2b_2b_3x_3+a_1b_2^2b_3^2-a_1b_2^2b_3x_3-a_1b_2b_3^2x_2+a_2^2a_3b_3x_1-b_1a_2b_2a_3b_3+b_1a_2b_2a_3x_3+b_1a_2b_2b_3^2-b_1a_2b_2b_3x_3-a_2b_2b_3^2x_1-a_2a_3f+b_2b_3f \right )\geq 0.
 \end{dmath}
\item The facet through points $v^1$, $v^4$, $v^7$ and $v_1^{11}$ is \ref{ineq17}.
\item The facet through points $v^4$, $v^5$, $v^7$ and $v_1^{11}$ is \ref{ineq16}.
\item The facet through points $v^1$, $v^7$, $v^8$ and $v_1^{11}$ is
\begin{dmath}
\label{ineq21b}
f-b_2b_3x_1-b_1a_3x_2-b_1a_2x_3+b_1a_2a_3+b_1b_2b_3 \geq 0.
\end{dmath}
\item The facet through points $v^5$, $v^7$, $v^8$ and $v_1^{11}$ is \ref{ineq2}.
\item The facet through points $v^1$, $v^3$, $v^4$ and $v_1^{11}$ is \ref{ineq1}.
\item The facet through points $v^3$, $v^4$, $v^5$ and $v_1^{11}$ is
\begin{dmath}
\label{ineq22b}
f-a_2a_3x_1-a_1b_3x_2-a_1b_2x_3+a_1a_2a_3+a_1b_2b_3\geq 0.
\end{dmath}
\end{itemize}

There are four inequalities that are not already contained in system $\T_1^{1-}$, we add these and in doing so obtain the system of inequalities that describes $\T_1^1 = \text{conv}(\Po_H \cup \{v_{1}^{11}\})$ (in case 1).

We now compute the additional volume of $\text{conv}(\T_1^1 \cup \{v_1^{12}\})$ compared to the volume of $\T_1^1$.  As we have done previously, we sum the volumes of $\text{conv}(\{v_1^{12}\} \cup F)$ for each facet, $F$, of $\T_1^1$ such that $v_1^{12}$ is beyond that facet.   We substitute $v_1^{12}$ into each relevant inequality (i.e., the system that describes $\T_1^1$) and if the result is negative then $v_1^{12}$ lies beyond that facet.  It is immediately clear that $v_1^{12}$ satisfies inequalities \ref{ineq1}, \ref{ineq2}, \ref{ineq7}, \ref{ineq9}-\ref{ineq18} and \ref{ineq21b}-\ref{ineq22b}.  We   also see immediately that inequalities \ref{ineq19b} and \ref{ineq20b} are violated.  To see that inequality \ref{ineq8} is also satisfied see appendix \S\ref{B:point11:ineq11}.

Therefore we see that $v_1^{12}$ is beyond two facets and we need to compute the volume of the convex hull of $v_1^{12}$ with each of these facets.

The extreme points that lie on the first facet are points $v^1, v^3, v^8$ and $v_1^{11}$.  The polytope $\text{conv}\{v^1, v^3, v^8, v_{1}^{11}, v_1^{12}\}$ is a 4-simplex with volume:
\begin{dmath*} a_2b_3(b_1-a_1)^2(b_2-a_2)^2(b_3-a_3)^2/\left(24(b_2b_3-a_2a_3)\right).\end{dmath*}

The extreme points that lie on the second facet are points $v^3, v^5, v^8$ and $v_1^{11}$.  The polytope $\text{conv}\{v^3, v^5, v^8, v_{1}^{11}, v_1^{12}\}$ is a 4-simplex with volume:
\begin{dmath*} a_2b_3(b_1-a_1)^2(b_2-a_2)^2(b_3-a_3)^2/\left(24(b_2b_3-a_2a_3)\right).\end{dmath*}

We now compute the additional facets; we take the four facets from adding each simplex and delete the facet that repeats.  This leaves us with the following six facet defining inequalities to compute:
\begin{itemize}
\item The facet through points $v^1$, $v^3$, $v^8$ and $v_1^{12}$ is \ref{ineq15}.
\item The facet through points $v^1$, $v^3$, $v_1^{11}$ and $v_1^{12}$ is \ref{ineq1}.
\item The facet through points $v^1$, $v^8$, $v_1^{11}$ and $v_1^{12}$ is \ref{ineq21b}.
\item The facet through points $v^3$, $v^5$, $v^8$ and $v_1^{12}$ is \ref{ineq18}.
\item The facet through points $v^3$, $v^5$, $v_1^{11}$ and $v_1^{12}$ is \ref{ineq22b}.
\item The facet through points $v^5$, $v^8$, $v_1^{11}$ and $v_1^{12}$ is \ref{ineq2}.
\end{itemize}

By adding and deleting the appropriate facets from system $H$ we see that we arrive at system 1.

Therefore, to compute the volume of $\Po_1$, we sum the volume of $\Po_H$ with that of the appropriate eight simplices, and we obtain our result for case 1.

\subsubsection[s2]{Case 2: $\bm{a_1b_2b_3-b_1a_2a_3 \leq 0}$}

In this case it is immediate to see that $v_1^{11}$ satisfies \ref{ineq3} and therefore lies beyond no further facets.  This means we now have a new polytope which is $\text{conv}(\Po_H \cup \{v_3^{11}\})$ (in case 2).  We refer to this polytope as $\T_1^2$, and we compute the facets of $\T_1^2$.

We begin with the facets of $\Po_H$ and delete the three facets that $v_1^{11}$ lies beyond (\ref{ineq4}-\ref{ineq6}).  Let us call this system $\T_1^{2-}$.  Now consider the four simplices we dealt with when computing the additional volume produced with $v_1^{11}$.  Each of these simplices has 5 facets; one of which corresponds to a deleted facet of $\Po_H$.\\

The remaining 4 facets of the first simplex are described by the planes through the following sets of points:
$\{v^1, v^4, v^5, v_1^{11}\}$, $\{v^1, v^4, v^7, v_1^{11}\}$, $\{v^1, v^5, v^7, v_1^{11}\}$ \\and $\{v^4, v^5, v^7, v_1^{11}\}$.

The remaining 4 facets of the second simplex are described by the planes through the following sets of points:
$\{v^1, v^5, v^7, v_1^{11}\}$, $\{v^1, v^5, v^8, v_1^{11}\}$, $\{v^1, v^7, v^8, v_1^{11}\}$ \\and $\{v^5, v^7, v^8, v_1^{11}\}$.

The remaining 4 facets of the third simplex are described by the planes through the following sets of points:
$\{v^1, v^3, v^4, v_1^{11}\}$, $\{v^1, v^3, v^5, v_1^{11}\}$, $\{v^1, v^4, v^5, v_1^{11}\}$ \\and $\{v^3, v^4, v^5, v_1^{11}\}$.\\

Consider these twelve facets and exclude the facets that are shared by more than one simplex.  This leaves eight facets.

 We   compute these eight facets to obtain the following:
\begin{itemize}
\item The facet through points $v^1$, $v^4$, $v^7$ and $v_1^{11}$ is \ref{ineq17}.
\item The facet through points $v^4$, $v^5$, $v^7$ and $v_1^{11}$ is \ref{ineq16}.
\item The facet through points $v^1$, $v^5$, $v^8$ and $v_1^{11}$ is
\begin{dmath}
\label{ineq23b}
\frac{1}{b_2(b_1-a_1)}\left (-a_1b_1a_2^2a_3+a_1b_1a_2a_3x_2+a_1b_1a_2b_2x_3-a_1b_1b_2^2b_3+a_1b_2^2b_3x_1+b_1a_2^2a_3x_1+b_1^2a_2b_2a_3-b_1^2a_2a_3x_2-b_1a_2b_2a_3x_1+b_1^2a_2b_2b_3-b_1^2a_2b_2x_3-b_1a_2b_2b_3x_1-a_1b_2f+b_1b_2f \right )\geq 0.
\end{dmath}
\item The facet through points $v^1$, $v^7$, $v^8$ and $v_1^{11}$ is \ref{ineq21b}.
\item The facet through points $v^5$, $v^7$, $v^8$ and $v_1^{11}$ is \ref{ineq2}.
\item The facet through points $v^1$, $v^3$, $v^4$ and $v_1^{11}$ is \ref{ineq1}.
\item The facet through points $v^1$, $v^3$, $v^5$ and $v_1^{11}$ is
 \begin{dmath}
\label{ineq24b}
\frac{1}{a_3(b_1-a_1)}\left (-a_1^2a_2a_3b_3-a_1^2b_2a_3b_3+a_1^2a_3b_3x_2+a_1^2b_2b_3x_3+a_1b_1a_2a_3^2+a_1a_2a_3b_3x_1-a_1b_1a_3b_3x_2+a_1b_2a_3b_3x_1+a_1b_1b_2b_3^2-a_1b_1b_2b_3x_3-a_1b_2b_3^2x_1-b_1a_2a_3^2x_1-a_1a_3f+b_1a_3f \right )\geq 0.
\end{dmath}
\item The facet through points $v^3$, $v^4$, $v^5$ and $v_1^{11}$ is \ref{ineq22b}.
\end{itemize}

There are four inequalities that are not already contained in system $\T_1^{2-}$; we add these and in doing this we obtain the system of inequalities that describes $T_1^2 = \text{conv}(\Po_H \cup \{v_1^{11}\})$ (in case 2).

We now compute the additional volume of $\text{conv}(\T_1^2 \cup \{v_1^{12}\})$ compared to the volume of $\T_1^2$.  As we have done previously, we sum the volumes  of $\text{conv}(\{v_1^{12}\} \cup F)$ for each facet, $F$, of $\T_1^2$ such that $v_1^{12}$ is beyond that facet.  We substitute $v_1^{12}$ into each relevant inequality (i.e., the system of inequalities that describes $T_1^2$ in case 2) and if the result is negative then $v_1^{12}$ lies beyond that facet.  It is immediately clear that $v_1^{12}$ satisfies inequalities \ref{ineq1}, \ref{ineq2}, \ref{ineq7}, \ref{ineq9}-\ref{ineq18}, \ref{ineq21b} and \ref{ineq22b}.  We   also see immediately that $v_1^{12}$ violates inequalities \ref{ineq3}, \ref{ineq23b} and \ref{ineq24b}.  To see that \ref{ineq8} is also satisfied see \S\ref{B:point11:ineq11}.

Therefore we know that $v_1^{12}$ is beyond three facets, and we need to compute the volume of the convex hull of $v_1^{12}$ with each of these facets.

The extreme points that lie on the first facet are points $v^1, v^3, v^5$ and $v^8$.  The polytope $\text{conv}\{v^1, v^3, v^5, v^{8}, v_1^{12}\}$ is a 4-simplex with volume:
\begin{dmath*}a_2b_3(b_1-a_1)^2(b_2-a_2)^2(b_3-a_3)^2/\left(24(b_2b_3-a_2a_3)\right).\end{dmath*}

The extreme points that lie on the second facet are points $v^1, v^5, v^8$ and $v_1^{11}$.  The polytope $\text{conv}\{v^1, v^5, v^8, v_1^{11}, v_1^{12}\}$ is a 4-simplex with volume:
\begin{dmath*} (b_1a_2(b_1-a_1)(b_2-a_2)^2(b_3-a_3)^3/\left(24(b_2b_3-a_2a_3)\right).\end{dmath*}

The extreme points that lie on the third and final facet are points $v^1, v^3, v^5$ and $v_1^{11}$.  The polytope $\text{conv}\{v^1, v^3, v^5, v_1^{11}, v_1^{12}\}$ is a 4-simplex with volume:
\begin{dmath*} a_1b_3(b_1-a_1)(b_2-a_2)^3(b_3-a_3)^2/\left(24(b_2b_3-a_2a_3)\right).\end{dmath*}

We now compute the additional facets; we take the four facets from adding each simplex and delete the three facets that are repeated.  This leaves us with the following six facet defining inequalities to compute:
\begin{itemize}
\item The facet through points $v^1$, $v^3$, $v^8$ and $v_1^{12}$ is \ref{ineq15}.
\item The facet through points $v^3$, $v^5$, $v^8$ and $v_1^{12}$ is \ref{ineq18}.
\item The facet through points $v^1$, $v^8$, $v_1^{11}$ and $v_1^{12}$ is \ref{ineq21b}.
\item The facet through points $v^5$, $v^8$, $v_1^{11}$ and $v_1^{12}$ is \ref{ineq2}.
\item The facet through points $v^1$, $v^3$, $v_1^{11}$ and $v_1^{12}$ is \ref{ineq1}.
\item The facet through points $v^3$, $v^5$, $v_1^{11}$ and $v_1^{12}$ is \ref{ineq22b}.
\end{itemize}

By adding and deleting the appropriate facets from system $H$ we see that we also arrive at system 1 in case 2.

Therefore, to compute the volume of $\Po_1$, we sum the volume of $\Po_H$ with that of the appropriate eight simplices, and we obtain our result for case 2.  \halmos

%% file: theorem_C.tex
\section{Proof of Thm.~\ref{TheoremC}.}
\label{linC}

A mapping from the proof of Theorem \ref{TheoremA} allows us to claim Theorem \ref{TheoremC} immediately. \halmos 

%% file: future.tex
\section{Possible extensions.} \label{sec:future}

Our results geometrically  quantify the tradeoff between different convexifications of
trilinear monomials.
Of course it would be nice to use our results to develop guidelines for attacking
trilinear monomials within an sBB code. In doing so, it should prove important to
develop guidelines for how our results could be applied to
formulations having many trilinear monomials overlapping on the same variables.
We have seen that our results are very robust
for scenarios where there is a high degree of overlap between trilinear monomials (see \cite{SpeakmanLee2016b}).
Also, we can imagine scoring each possible relaxation according to its volume, and then
aggregating the scores to decide on what to do for each trilinear monomial.
Another important issue is how to effectively make branching decisions in the context of our relaxations.
Guided by our volume results, we have made some significant progress in this direction
(see \cite{SpeakmanLee2016a}).

It would be natural and certainly difficult
 to extend our work to multilinear monomials having $n>3$.
In particular, advances for the
important case of $n=4$ could have immediate impact;
\cite{CafieriLeeLiberti10} found, via experiments, that composing
a trilinear and bilinear convexification in the
manner suggested by $(x_i x_j)x_k x_l$ was a
good strategy. They further observed sensitivity to
the bounds on the variables, but they reached no clear conclusion
on how to factor in that aspect.
Restricting to this type of convexification,
we could apply our results by substituting $w\in[a_i a_j, b_i b_j]$
to arrive at the trilinear monomial $w x_k x_l$, which can then be
analyzed and relaxed according to our methodology.
Of course, for a general quadrilinear monomial, there are six
choices of which pair of variables will be treated as $\{x_i,x_j\}$,
so we can analyze all six possibilities and take the best overall.

Also, there is the possibility of extending our results
on trilinear monomials to
(i) box domains that are not necessarily nonnegative,
(ii) domains other than boxes, and
(iii) other low-dimensional functions.

We hope that our work is just a first step in using
volume to better understand and mathematically quantify the tradeoffs
involved in developing sBB strategies for factorable formulations.


%% file: Appendix.tex
\section{Appendix.}\label{app}



Throughout the proofs, we have repeatedly claimed that certain quantities are nonnegative for any choice of $a_1,a_2,a_3,b_1,b_2,b_3$, such that, $0<a_i<b_i$, for all $i$ and
\[
a_1b_2b_3+b_1a_2a_3 \leq b_1a_2b_3 + a_1b_2a_3 \leq b_1b_2a_3 + a_1a_2b_3.
\]
In this appendix, we provide proofs for the cases that are not immediate.
As will become apparent, we need to demonstrate that many different 6-variable polynomials are nonnegative on the relevant parameter space. Generally, such demonstrations can be tricky global-optimization problems, and in many cases  sum-of-squares proofs are not available; rather, we often make somewhat ad hoc arguments. Still, we can place some efficiency on all of this by establishing some technical lemmas.

\subsection[9.1]{}
\label{lemmas}
We begin with the following lemmas that will be helpful in establishing the nonnegativity of certain quantities:\\

\begin{lemma}\label{lem0}
For all choices of parameters that meet our assumptions we have: $b_1a_2-a_1b_2\geq 0$, $b_1a_3-a_1b_3 \geq 0$ and $b_2a_3-a_2b_3 \geq 0$.
\end{lemma}

\proof{Proof.}
$(b_3-a_3)(b_1a_2-a_1b_2) = b_1a_2b_3 + a_1b_2a_3 - a_1b_2b_3 -b_1a_2a_3 \geq 0$
 by our original assumptions \ref{Omega}.
This implies $b_1a_2-a_1b_2 \geq 0$, because $b_3-a_3 > 0$.
$b_1a_3-a_1b_3 \geq 0$ and $b_2a_3-a_2b_3 \geq 0$ follow from \ref{Omega} in a similar way. \halmos \endproof

\begin{lemma}\label{lem1}
Let $A,B,C,D,E,F \in \R$ with $A \geq B \geq C \geq 0$ and $D \geq 0$, $E\geq 0$, $F\leq0$. Also let, $D+E+F=0$.  Then $AD+BE+CF \geq 0$.
\end{lemma}

\proof{Proof.}
$AD+BE+CF = AD + BE - C(D+E) \geq BD + BE - CD - CE = (B-C)(D+E) \geq 0$.
\halmos \endproof

\begin{lemma}
\label{lem2}
Let $A,B,C,D,E,F \in \R$ with $A \geq B \geq C \geq 0$ and $D \geq 0$, $E\leq 0$, $F\leq0$. Also let, $D+E+F=0$.  Then $AD+BE+CF \geq 0$.
\end{lemma}

\proof{Proof.}
$AD+BE+CF = -A(E+F) + BE + CF = E(B-A) + F(C-A) \geq 0$.
\halmos \endproof

\begin{lemma}
\label{lem3}
Let $A,B,C,D \in \R$ with $A \geq B\geq 0$, $C+D \geq 0$, $C \geq 0$.  Then $AC+BD \geq0$.
\end{lemma}

\proof{Proof.}
$AC+BD \geq B(C+D) \geq 0$.
\halmos \endproof

\subsection[9.2]{}
\label{A:midpoint911:ineq9}
Substituting $M(v_3^{9},v_3^{11})$ into inequality \ref{ineq9} of the convex hull we obtain:
\begin{dmath*}
\begin{aligned}
&\left( -a_1^2a_2^2a_3+a_1^2a_2b_2b_3+a_1a_2^2a_3b_1+2a_1a_2a_3b_1b_2-2a_1a_2b_1b_2b_3\right.\\
&\qquad\left.-a_1a_3b_1b_2^2-a_2^2a_3b_1^2+a_2^2b_1^2b_3-a_2b_1^2b_2b_3+b_1^2b_2^2b_3\right ) \Big/2(b_1b_2-a_1a_2),
\end{aligned}
\end{dmath*}
the numerator of which can be rewritten as
\[
b_1b_2\Big ((b_1b_3-a_1a_3)(b_2-a_2)\Big)+b_1a_2\Big((b_1a_2-a_1b_2)(b_3-a_3)\Big)+a_1a_2\Big((b_2b_3-a_2a_3)(a_1-b_1)\Big),
\]
and is nonnegative by Lemmas \ref{lem0} and \ref{lem1}.

\subsection[9.3]{}
\label{A:midpoint1012:ineq10}
Substituting $M(v_3^{10},v_3^{12})$ into inequality \ref{ineq10} of the convex hull we obtain:
\begin{equation*}
\begin{aligned}
&\left(-a_1^2a_2^2a_3+a_1^2a_2a_3b_2-a_1^2a_3b_2^2+a_1^2b_2^2b_3+a_1a_2^2b_1b_3+2a_1a_2a_3b_1b_2\right.\\
&\qquad\left.-2a_1a_2b_1b_2b_3-a_1b_1b_2^2b_3-a_2a_3b_1^2b_2+b_1^2b_2^2b_3\right ) \Big/ 2(b_1b_2-a_1a_2),
\end{aligned}
\end{equation*}
the numerator of which can be rewritten as
\begin{equation*}
b_1b_2\Big((b_2b_3-a_2a_3)(b_1-a_1)\Big)+a_1b_2\Big((b_1a_2-a_1b_2)(a_3-b_3)\Big)+a_1a_2\Big((b_1b_3-a_1a_3)(a_2-b_2)\Big),
\end{equation*}
and is nonnegative by Lemmas \ref{lem0} and \ref{lem2}.

\subsection[9.4]{}
\label{A:point11:ineq9}
Substituting point $v_3^{11}$ into inequality \ref{ineq9} of the convex hull or substituting $v_3^{12}$ into inequality \ref{ineq10} we obtain:
\begin{dmath*}
\begin{aligned}
&\left(-a_1^2a_2^2a_3+a_1^2a_2b_2b_3+a_1a_2^2b_1b_3+3a_1a_2a_3b_1b_2-3a_1a_2b_1b_2b_3\right.\\
&\qquad\left.-a_1a_3b_1b_2^2-a_2a_3b_1^2b_2+b_1^2b_2^2b_3\right) \Big/ (b_1b_2-a_1a_2).
\end{aligned}
\end{dmath*}
The numerator can be rewritten as
\begin{dmath*}
b_3\Big(a_1a_2(a_1b_2-b_1b_2)+b_1a_2(a_1a_2-a_1b_2)+b_1b_2(b_1b_2-a_1a_2)\Big)+a_3\Big(a_1a_2(b_1b_2-a_1a_2)+b_1a_2(a_1b_2-b_1b_2)+b_1b_2(a_1a_2-a_1b_2)\Big) \equalscolon  b_3Y + a_3Z.
\end{dmath*}
Then we can see $Y + Z = (b_2-a_2)(b_1-a_1)(b_1b_2-a_1a_2)$, which is nonnegative.
Furthermore, by Lemma \ref{lem2} we have $Y \geq 0$. Therefore, by Lemma \ref{lem3} the numerator is nonnegative.


\subsection[9.5]{}
\label{B:midpoint910:ineq3}
Substituting   $M(v_1^{9},v_1^{10})$ into inequality \ref{ineq3} of the convex hull we obtain:
\begin{dmath*}
\begin{aligned}
&\left(-2a_1a_2^2a_3b_3+a_1a_2^2b_3^2-a_1a_2a_3^2b_2+4a_1a_2a_3b_2b_3-a_1a_2b_2b_3^2-a_1b_2^2b_3^2+a_2^2a_3^2b_1\right.\\
&\qquad\left. +a_2^2a_3b_1b_3-a_2^2b_1b_3^2-4a_2a_3b_1b_2b_3+2a_2b_1b_2b_3^2+a_3b_1b_2^2b_3\right)
\Big/ 2(b_2b_3-a_2a_3),
\end{aligned}
\end{dmath*}
the numerator of which can be rewritten as
\begin{dmath*}
b_2b_3\Big(b_2(b_1a_3-a_1b_3)+a_2(b_1b_3-2b_1a_3+a_1a_3) \Big)+ a_2b_3\Big((b_1-a_1)(b_2-a_2)(b_3-a_3) \Big) + a_2a_3\Big(b_2(2a_1b_3-b_1b_3-a_1a_3) +a_2(b_1a_3-a_1b_3) \Big) \equalscolon b_2b_3X + a_2b_3Y + a_2a_3Z.
\end{dmath*}
Now, we write $X \equalscolon b_2V+a_2W$ and see that $V+W=(b_1-a_1)(b_3-a_3) \geq 0$ and, by Lemma \ref{lem0}, $V\geq 0$.  Therefore $X\geq0$ by Lemma \ref{lem3}.  Because $X+Y+Z=0$, by Lemma \ref{lem1} we have that the numerator is nonnegative.


\subsection[9.6]{}
\label{B:midpoint910:ineq4}
Substituting $M(v_1^{9},v_1^{10})$ into inequality \ref{ineq4} of the convex hull, we obtain:
\begin{dmath*}
 \begin{aligned}
& \left(-a_1a_2^2a_3b_3-2a_1a_2a_3^2b_2+4a_1a_2a_3b_2b_3+a_1a_3^2b_2^2-a_1a_3b_2^2b_3-a_1b_2^2b_3^2+a_2^2a_3^2b_1\right.\\
&\qquad \left. +a_2a_3^2b_1b_2-4a_2a_3b_1b_2b_3+a_2b_1b_2b_3^2-a_3^2b_1b_2^2+2a_3b_1b_2^2b_3\right)\Big/2(b_2b_3-a_2a_3),
\end{aligned}
\end{dmath*}
the numerator of which can be rewritten as
\begin{dmath*}
 b_2b_3\Big(b_3(b_1a_2-a_1b_2)+a_3(b_1b_2-2b_1a_2+a_1a_2) \Big)+ b_2a_3\Big((b_1-a_1)(b_2-a_2)(b_3-a_3) \Big) + a_2a_3\Big(b_3(2a_1b_2-b_1b_2-a_1a_2) +a_3(b_1a_2-a_1b_2) \Big) \equalscolon b_2b_3X + b_2a_3Y + a_2a_3Z.
 \end{dmath*}
 Now, we write $X \equalscolon b_3V+a_3W$ and see that $V+W=(b_1-a_1)(b_2-a_2) \geq 0$ and, by Lemma \ref{lem0} $V\geq 0$.  Therefore $X\geq0$ by Lemma \ref{lem3}.  Because $X+Y+Z=0$, by Lemma \ref{lem1} we have that the numerator is nonnegative.
%

\subsection[9.7]{}
\label{B:midpoint911:ineq4}
Substituting $M(v_1^{9},v_1^{11})$ into inequality \ref{ineq4} of the convex hull, we obtain:
\begin{dmath*}
 \begin{aligned}
& \left(-a_1a_2^2a_3b_3+2a_1a_2a_3b_2b_3-a_1a_3^2b_2^2+a_1a_3b_2^2b_3-a_1b_2^2b_3^2+a_2^2a_3^2b_1-a_2a_3^2b_1b_2\right.\\
&\qquad \left.-2a_2a_3b_1b_2b_3+a_2b_1b_2b_3^2+a_3^2b_1b_2^2 \right) \Big/2(b_2b_3-a_2a_3),
 \end{aligned}
\end{dmath*}
the numerator of which can be rewritten as
\begin{equation*}
b_2b_3\Big((b_1a_2-a_1b_2)(b_3-a_3)\Big) + b_2a_3\Big((b_2a_3-a_2b_3)(b_1-a_1)\Big) + a_2a_3\Big((b_1a_3-a_1b_3)(a_2-b_2)\Big),
\end{equation*}
 which is nonnegative by Lemmas \ref{lem1} and \ref{lem0}.

\subsection[9.8]{}
\label{B:midpoint911:ineq5}
Substituting $M(v_1^{9},v_1^{11})$ into inequality \ref{ineq5} of the convex hull we obtain:
\begin{dmath*}
\left(b_2-a_2\right)\left(a_1a_2a_3b_3-a_1b_2b_3^2-a_2a_3^2b_1-a_2a_3b_1b_3+a_2b_1b_3^2+a_3^2b_1b_2 \right)
\Big/2(b_2b_3-a_2a_3),
\end{dmath*}
where the second multiplicand of the numerator can be rewritten as
\begin{dmath*}b_3\Big (b_3(b_1a_2-a_1b_2) \Big) + a_3 \Big(a_1a_2b_3-b_1a_2a_3+b_1b_2a_3-b_1a_2b_3 \Big) \equalscolon b_3Y + a_3Z,\end{dmath*}
now $Y+Z = (b_1a_3-a_1b_3)(b_2-a_2) \geq 0$ (Lemma \ref{lem0}), and $Y \geq 0$ (Lemma \ref{lem0}), therefore by Lemma \ref{lem3} we have that $b_3Y + a_3Z$ is nonnegative.

\subsection[9.9]{}
\label{B:midpoint911:ineq6}
Substituting $M(v_1^{9},v_1^{11})$ into inequality \ref{ineq6} of the convex hull we obtain:
\begin{dmath*}
\left(b_3-a_3 \right) \left( a_1a_2^2b_3-a_1a_2b_2b_3+a_1a_3b_2^2-a_1b_2^2b_3-a_2^2a_3b_1+a_2b_1b_2b_3\right)\Big/2(b_2b_3-a_2a_3),
\end{dmath*}
where the second multiplicand of the numerator can be rewritten as
\begin{dmath*}b_2\Big(b_3(b_1a_2-a_1b_2)+a_1(b_2a_3-a_2b_3) \Big) + a_2\Big(a_2(a_1b_3-b_1a_3) \Big)\equalscolon b_2Y + a_2Z,\end{dmath*}
where $Y+Z = (b_1a_2-a_1b_2)(b_3-a_3) \geq 0$ and $Y \geq 0$ (both Lemma \ref{lem0}), therefore by Lemma \ref{lem3} we have that this term is nonnegative.

\subsection[9.10]{}
\label{B:midpoint911:ineq11}
Substituting $M(v_1^{9},v_1^{11})$ into inequality \ref{ineq11} of the convex hull we obtain:
\begin{dmath*}
\begin{aligned}
\left(-a_1a_2^2a_3^2+a_1a_2^2a_3b_3-a_1a_2^2b_3^2+2a_1a_2a_3b_2b_3-a_1a_3b_2^2b_3+a_2^2b_1b_3^2\right.\\
\qquad \left.+a_2a_3^2b_1b_2-2a_2a_3b_1b_2b_3-a_2b_1b_2b_3^2+b_1b_2^2b_3^2\right) \Big / 2(b_2b_3-a_2a_3),
\end{aligned}
\end{dmath*}
the numerator of which can be rewritten as
\begin{equation*}b_2b_3\Big((b_1b_3-a_1a_3)(b_2-a_2)\Big) + a_2b_3\Big((b_2a_3-a_2b_3)(a_1-b_1)\Big) + a_2a_3\Big((b_1b_2-a_1a_2)(a_3-b_3)\Big),\end{equation*}
 which is nonnegative by Lemmas \ref{lem2} and \ref{lem0}.

\subsection[9.11]{}
\label{B:midpoint1012:ineq3}
Substituting $M(v_1^{10},v_1^{12})$ into inequality \ref{ineq3} of the convex hull we obtain:
\begin{dmath*}
\begin{aligned}
\left(-a_1a_2^2b_3^2-a_1a_2a_3^2b_2+2a_1a_2a_3b_2b_3+a_1a_2b_2b_3^2-a_1b_2^2b_3^2+a_2^2a_3^2b_1\right.\\
\qquad \left.-a_2^2a_3b_1b_3+a_2^2b_1b_3^2-2a_2a_3b_1b_2b_3+a_3b_1b_2^2b_3\right)  \Big / 2(b_2b_3-a_2a_3),
\end{aligned}
\end{dmath*}
the numerator of which can be rewritten as
\begin{equation*}b_2b_3\Big((b_1a_3-a_1b_3)(b_2-a_2)\Big) + a_2b_3\Big((b_2a_3-a_2b_3)(a_1-b_1)\Big) + a_2a_3\Big((b_1a_2-a_1b_2)(a_3-b_3)\Big),\end{equation*}
 which is nonnegative by Lemma \ref{lem2} and Lemma \ref{lem0}.

\subsection[9.12]{}
\label{B:midpoint1012:ineq5}
Substituting $M(v_1^{10},v_1^{12})$ into inequality \ref{ineq5} of the convex hull we obtain:
\begin{dmath*}
\left(b_3-a_3 \right) \left( a_1a_2a_3b_2-a_1b_2^2b_3-a_2^2a_3b_1+a_2^2b_1b_3-a_2a_3b_1b_2+a_3b_1b_2^2\right) \Big / 2(b_2b_3-a_2a_3),
\end{dmath*}
where the second multiplicand of the numerator can be rewritten as
\begin{dmath*} b_2\Big( b_2(b_1a_3-a_1b_3)\Big) + a_2\Big(b_1(a_2b_3-b_2a_3)+a_3(a_1b_2-b_1a_2)\Big)\equalscolon b_2Y + a_2Z,\end{dmath*}
where $Y+Z = (b_1a_2-a_1b_2)(b_3-a_3) \geq 0$ and $Y \geq 0$ (both by Lemma \ref{lem0}).  Therefore by Lemma \ref{lem3} we have that $b_2Y + a_2Z$ is nonnegative.

\subsection[9.13]{}
\label{B:midpoint1012:ineq6}
Substituting $M(v_1^{10},v_1^{12})$ into inequality \ref{ineq6} of the convex hull we obtain:
\begin{dmath*}
\left(b_2-a_2 \right) \left(a_1a_2b_3^2+a_1a_3^2b_2-a_1a_3b_2b_3-a_1b_2b_3^2-a_2a_3^2b_1+a_3b_1b_2b_3 \right) \Big / 2(b_2b_3-a_2a_3),
\end{dmath*}
where the second multiplicand of the numerator can be rewritten as
\begin{dmath*} b_3\Big( a_1(a_2b_3-b_2a_3)+b_2(b_1a_3-a_1b_3)\Big) + a_3\Big(a_3(a_1b_2-b_1a_2)\Big)\equalscolon b_3Y + a_3Z,\end{dmath*}
where $Y+Z = (b_1a_3-a_1b_3)(b_2-a_2) \geq 0$ and $Z \leq 0 \Longrightarrow Y \geq 0$ (both by Lemma \ref{lem0}).  Therefore by Lemma \ref{lem3} we have that $b_3Y + a_3Z$ is nonnegative.

\subsection[9.14]{}
\label{B:midpoint1012:ineq8}
Substituting $M(v_1^{10},v_1^{12})$ into inequality \ref{ineq8} of the convex hull we obtain:
\begin{dmath*}
\begin{aligned}
\left(-a_1a_2^2a_3^2+a_1a_2a_3^2b_2+2a_1a_2a_3b_2b_3-a_1a_2b_2b_3^2-a_1a_3^2b_2^2+a_2^2a_3b_1b_3\right.\\
\left.-2a_2a_3b_1b_2b_3+a_3^2b_1b_2^2-a_3b_1b_2^2b_3+b_1b_2^2b_3^2\right) \Big / 2(b_2b_3-a_2a_3),
\end{aligned}
\end{dmath*}
the numerator of which simplifies to
\begin{equation*}
b_2b_3\Big((b_1b_2-a_1a_2)(b_3-a_3)\Big) + b_2a_3\Big((b_2a_3-a_2b_3)(b_1-a_1)\Big) + a_2a_3\Big((b_1b_3-a_1a_3)(a_2-b_2)\Big),
\end{equation*}
and is nonnegative by Lemmas \ref{lem0} and \ref{lem1}.

\subsection[9.15]{}
\label{B:point9:ineq4}
Substituting point $v_1^{9}$ into inequality \ref{ineq4} of the convex hull or substituting point $v_1^{10}$ into inequality \ref{ineq3} we obtain:
\begin{dmath*}
\begin{aligned}
\left (-a_1a_2^2a_3b_3-a_1a_2a_3^2b_2+3a_1a_2a_3b_2b_3-a_1b_2^2b_3^2+a_2^2a_3^2b_1\right.\\
\left.-3a_2a_3b_1b_2b_3+a_2b_1b_2b_3^2+a_3b_1b_2^2b_3\right ) \Big / (b_2b_3-a_2a_3),
\end{aligned}
\end{dmath*}
the numerator of which can be rewritten as
\begin{dmath*} b_3\Big (b_2\big(b_2(b_1a_3-a_1b_3) + a_2(a_1a_3+b_1b_3-2b_1a_3) \big)\Big) + a_3\Big (a_2 \big(b_2(2a_1b_3-a_1a_3-b_1b_3)+a_2(b_1a_3-a_1b_3) \big)\Big)\equalscolon b_3Y + a_3Z.\end{dmath*}
Now, we write $Y \equalscolon b_2^2V+a_2b_2W$ and see that $V+W=(b_1-a_1)(b_3-a_3) \geq 0$ and, by Lemma \ref{lem0} $V\geq 0$.  Therefore $Y\geq 0$ by Lemma \ref{lem3}.  Because $Y+Z = (b_1a_3-a_1b_3)(b_2-a_2)^2 \geq 0$ (Lemma \ref{lem0}), by Lemma \ref{lem3} we have that the numerator is nonnegative.

\subsection[9.16]{}
\label{B:point11:ineq11}
 Substituting point $v_1^{11}$ into inequality \ref{ineq11} of the convex hull or substituting point $v_1^{12}$ into \ref{ineq8} we obtain:
 \begin{dmath*}
\begin{aligned}
\left(-a_1a_2^2a_3^2+3a_1a_2a_3b_2b_3-a_1a_2b_2b_3^2-a_1a_3b_2^2b_3+a_2^2a_3b_1b_3\right.\\
\left.+a_2a_3^2b_1b_2-3a_2a_3b_1b_2b_3+b_1b_2^2b_3^2\right ) \Big /  (b_2b_3-a_2a_3),
\end{aligned}
\end{dmath*}
the numerator of which can be rewritten as
\begin{dmath*}  b_1\Big(b_2b_3(b_2b_3-a_2a_3) + b_2a_3(a_2a_3-a_2b_3) + a_2a_3(a_2b_3-b_2b_3)\Big) + a_1\Big(b_2b_3(a_2a_3-a_2b_3) + b_2a_3(a_2b_3-b_2b_3) + a_2a_3(b_2b_3-a_2a_3)\Big)\equalscolon b_1Y + a_1Z,\end{dmath*}
where: $Y + Z = (b_2b_3-a_2a_3)(b_2-a_2)(b_3-a_3) \geq 0$, and by Lemma \ref{lem2} we have $Y\geq 0$.  Therefore by Lemma \ref{lem3} we have that the numerator is nonnegative.

\subsection[9.17]{}
\label{B:midpoint910:ineq5}
 Substituting $M(v_1^{9},v_1^{10})$ into inequality \ref{ineq5} of the convex hull we obtain
 \begin{dmath*}
 (2b_2b_3-a_2b_3-b_2a_3)(a_1a_2a_3-a_1b_2b_3-2b_1a_2a_3+b_1a_2b_3+b_1b_2a_3)\Big / 2(b_2b_3-a_2a_3),
 \end{dmath*}
where the second multiplicand of the numerator can be rewritten as
\begin{dmath*}
b_2(b_1a_3-a_1b_3)+a_2(a_1a_3-2b_1a_3+b_1b_3),
\end{dmath*}
which is nonnegative by Lemmas \ref{lem0} and \ref{lem3}.

\subsection[9.18]{}
\label{B:midpoint910:ineq6}
Substituting $M(v_1^{9},v_1^{10})$ into inequality \ref{ineq6} of the convex hull we obtain
 \begin{dmath*}
 (a_2b_3+b_2a_3-2a_2a_3)(a_1a_2b_3+a_1b_2a_3-2a_1b_2b_3-b_1a_2a_3+b_1b_2b_3)\Big / 2(b_2b_3-a_2a_3),
 \end{dmath*}
where the second multiplicand of the numerator can be rewritten as
\begin{dmath*}
b_2(b_1b_3-2a_1b_3+a_1a_3)+a_2(a_1b_3-b_1a_3),
\end{dmath*}
which is nonnegative by Lemmas \ref{lem0} and \ref{lem3}.

\subsection[9.19]{}
\label{B:point9:ineq5}
Substituting point $v_1^{9}$ into inequality \ref{ineq5} of the convex hull we obtain
\begin{dmath*}
b_3(b_2-a_2)\Big(b_2(b_1a_3-a_1b_3)+a_2(a_1a_3-2b_1a_3+b_1b_3)\Big) \Big / (b_2b_3-a_2a_3),
\end{dmath*}
which is nonnegative by Lemmas \ref{lem0} and \ref{lem3}.

\subsection[9.20]{}
\label{B:point9:ineq6}
Substituting point $v_1^{9}$ into inequality \ref{ineq6} of the convex hull we obtain
\begin{dmath*}
a_2(b_3-a_3)\Big(b_2(b_1b_3-2a_1b_3+a_1a_3)+a_2(a_1b_3-b_1a_3)\Big) \Big / (b_2b_3-a_2a_3),
\end{dmath*}
which is nonnegative by Lemmas \ref{lem0} and \ref{lem3}.

\subsection[9.21]{}
\label{B:point10:ineq5}
Substituting point $v_1^{10}$ into inequality \ref{ineq5} of the convex hull we obtain
\begin{dmath*}
b_2(b_3-a_3)\Big(b_2(b_1a_3-a_1b_3)+a_2(a_1a_3-2b_1a_3+b_1b_3)\Big) \Big / (b_2b_3-a_2a_3),
\end{dmath*}
which is nonnegative by Lemmas \ref{lem0} and \ref{lem3}.

\subsection[9.22]{}
\label{B:point10:ineq6}
Substituting point $v_1^{10}$ into inequality \ref{ineq6} of the convex hull we obtain
\begin{dmath*}
a_3(b_2-a_2)\Big(b_2(b_1b_3-2a_1b_3+a_1a_3)+a_2(a_1b_3-b_1a_3)\Big) \Big / (b_2b_3-a_2a_3),
\end{dmath*}
which is nonnegative by Lemmas \ref{lem0} and \ref{lem3}.

%
%
%
%
%
